\documentclass[a4paper,11pt,reqno]{amsart}
\usepackage{graphicx, psfrag, amsmath, amscd, amssymb, array}
\newtheorem{lem}{Lemma}[section]
\newtheorem{thm}[lem]{Theorem}
\newtheorem{pro}[lem]{Proposition}
\newtheorem{cor}[lem]{Corollary}
\newtheorem{exa}[lem]{Example}
\newtheorem{con}[lem]{Conjecture}
\newtheorem{defi}[lem]{Definition}

\numberwithin{equation}{section}

\newcommand{\ZZ}{{\mathbb{Z}}}

\newcommand{\sign}{{\textsf{sign}}}

\newcommand{\F}{{\mathcal{F}}}

\newcommand{\SSS}{{\mathcal{S}}}
\newcommand{\RS}{{\mathcal{RS}}}
\newcommand{\T}{{\mathcal{T}}}

\newcommand{\rmlab}{{\textsf{rmlab}}}

\newcommand{\emp}{{\textsf{emp}}}

\newcommand{\npk}{{\textsf{npk}}}

\newcommand{\nva}{{\textsf{nva}}}
\newcommand{\w}{{\textsf{w}}}

\newcommand{\gae}{{\textsf{gae}}}

\newcommand{\AD}{{\mathcal{AD}}}

\newcommand{\I}{{\textrm{I}}}
\newcommand{\II}{{\textrm{II}}}

\title[Arnold Families Revisited]{Springer Numbers and Arnold Families Revisited}

\author[{S.-P. Eu}]{Sen-Peng Eu}
\address{Department of Mathematics, National Taiwan Normal University, Taipei 116325, and Chinese Air Force Academy, Kaohsiung 820009, Taiwan, ROC}
\email{speu@math.ntnu.edu.tw}

\author[T.-S. Fu]{Tung-Shan Fu}
\address{Department of Applied Mathematics, National Pingtung University, Pingtung 900391, Taiwan, ROC}
\email{tsfu@mail.nptu.edu.tw}

\begin{document}


\begin{abstract} For the calculation of Springer numbers (of root systems) of type $B_n$ and $D_n$, Arnold introduced a signed analogue of alternating permutations, called $\beta_n$-snakes, and derived recurrence relations for enumerating the $\beta_n$-snakes starting with $k$. The  results are presented in the form of double triangular arrays ($v_{n,k}$) of integers, $1\le |k|\le n$. An Arnold family is a sequence of sets of such objects as $\beta_n$-snakes that are counted by $(v_{n,k})$. As a refinement of Arnold's result, we give analogous arrays of polynomials, defined by recurrence, for the calculation of the polynomials associated with successive derivatives of $\tan x$ and $\sec x$, established by Hoffman. Moreover, we provide some new Arnold families of combinatorial objects that realize the polynomial arrays, which are signed variants of Andr\'{e} permutations and Simsun permutations.
\end{abstract}

\maketitle

\section{Introduction} 

Let $\mathfrak{S}_n$ be the symmetric group on $[n]:=\{1,2,\ldots,n\}$.
The classical \emph{Euler numbers} $E_n$, defined by
\[
1+\sum_{n\ge 1}E_n\frac{x^n}{n!}=\tan x+\sec x,
\]
count the \emph{alternating permutations} in $\mathfrak{S}_n$, i.e., $\sigma\in\mathfrak{S}_n$ such that
$\sigma_1>\sigma_2<\sigma_3>\cdots \sigma_n$ \cite[A000111]{oeis}. Among other objects, the Euler numbers also count \emph{Andr\'{e} permutations of the first and second kinds} \cite{FS} and \emph{Simsun permutations}, which are in connection with the $cd$-index of the symmetric group (see \cite{Stanley}). In 1877, Seidel defined the triangular array ($E_{n,k}$) for the calculation of $E_n$ by the recurrence
\begin{equation}
E_{n,k}=E_{n,k-1}+E_{n-1,n-k+1}\quad (n\ge k\ge 2)
\end{equation}
with $E_{1,1}=1$, $E_{n,1}=0$ ($n\ge 2$), and proved that $E_n=E_{n,1}+E_{n,2}+\cdots+E_{n,n}$. Entringer \cite{Entringer} proved that $E_{n,k}$ is the number of alternating permutations in $\mathfrak{S}_n$ starting with $k$. 

As the symmetric group is the Coxeter group of type $A_n$, Springer defined an integer $K(W)$ algebraically in a unified manner for each Coxeter group $W$ such that $K(\mathfrak{S}_n)=E_n$ (e.g., \cite{JNT, Springer}) 
For the groups of signed permutations of type $B_n$ and $D_n$, 
Arnold \cite{Arnold-1} introduced a signed analogue of alternating permutations, called $\beta_n$-\emph{snakes}, counted by the numbers $K(B_n)$ and $K(D_n)$. Several of the initial values are listed in Table \ref{tab:Euler-Springer-numbers}.
\begin{table}[ht]
\caption{The initial values of $K(W)$ for $W=\mathfrak{S}_n, B_n$ and $D_n$.}
\begin{tabular}{c|cccccccc}   \hline
$n$
   &  1 & 2 & 3 & 4 &  5 &  6  & 7 & 8 \\
    \hline
   $K(\mathfrak{S}_n)$     & 1 & 1 & 2  &  5  &  16   &  61   & 272   & 1385  \\
   $K(B_n)$  & 1 & 3 & 11 &  57 &  361  &  2763 & 24611 & 250737\\
   $K(D_n)$  & 1 & 1 & 5  &  23 &  151  &  1141 & 10205 & 103823\\
      \hline    
\end{tabular}
\label{tab:Euler-Springer-numbers}
\end{table}

Arnold derived recurrence relations for the enumeration of $\beta_n$-snakes with respect to the first entry, presenting these integers in triangular arrays like Seidel--Entringer triangle for $E_n$ (e.g., \cite{GSZ, Stanley}). Namely, he defined the following double triangular arrays ($v_{n,k}$), $1\le |k|\le n$, by the recurrence relations
\begin{equation} \label{eqn:Springer-recurrence}
\begin{aligned}
v_{n,-k} &= v_{n,-k-1}+v_{n-1,k} \quad (n>k\ge 1),\\
v_{n,1}  &= v_{n,-1} \quad (n\ge 2),\\
v_{n,k}  &= v_{n,k-1}+v_{n-1,-k+1} \quad (n\ge k>1),
\end{aligned}
\end{equation}
with $v_{1,1}=v_{1,-1}=1$, $v_{n,-n}=0$ ($n\ge 2$), and proved that $K(B_n)=v_{n,1}+v_{n,2}+\cdots+v_{n,n}$ and $K(D_n)=v_{n,-1}+v_{n,-2}+\cdots+v_{n,-n}$ (see Theorem \ref{thm:Arnold-Theorem}). Several of the initial values of $v_{n,k}$ are listed in Table \ref{tab:Arnold-Springer-numbers}. These integers $K(B_n)$ and $K(D_n)$ are known as \emph{Springer numbers of type} $B$ and $D$, respectively.
% for $1\le |k|\le n$, where
% \[
% \SSS_{n,k}:=\{\sigma\in\SSS_n : \sigma_1=k\}.
% \]

{\small
\[
\begin{array}{ccc}
\begin{array}{ccccccc}
              &   &    & v_{1,-1} &   &   &    \\
              &   & v_{2,2}   & \leftarrow   &  v_{2,1}  &    &    \\
              & v_{3,-3} & \rightarrow & v_{3,-2}   & \rightarrow   & v_{3,-1}   &   \\
     v_{4,4}  & \leftarrow &  v_{4,3} &  \leftarrow  & v_{4,2}  & \leftarrow &  v_{4,1}  \\
              &   &    & \cdots &   &   &    \\          
\end{array}
& &
\begin{array}{ccccccc}
              &   &    & v_{1,1} &   &   &    \\
              &   & v_{2,-1}   & \leftarrow   &  v_{2,-2}  &    &    \\
              & v_{3,1} & \rightarrow & v_{3,2}   & \rightarrow   & v_{3,3}   &   \\
     v_{4,-1}  & \leftarrow &  v_{4,-2} &  \leftarrow  & v_{4,-3}  & \leftarrow &  v_{4,-4}  \\
              &   &    & \cdots &   &   &    \\            
\end{array} 
            \\[2ex]
& & \\
\Updownarrow & & \Updownarrow \\
& & \\
\begin{array}{ccccccc}
              &   &    & 1 &   &   &    \\
              &   & 2   & \leftarrow   &  1  &    &    \\
              & 0 & \rightarrow & 2   & \rightarrow   & 3   &   \\
     16  & \leftarrow &  16 &  \leftarrow  & 14  & \leftarrow &  11 \\
              &   &    & \cdots &   &   &       
\end{array}
& &
\begin{array}{ccccccc}
              &   &    & 1 &   &   &    \\
              &   & 1   & \leftarrow   &  0  &    &    \\
              & 3 & \rightarrow & 4   & \rightarrow   & 4   &   \\
     11  & \leftarrow &  8 &  \leftarrow  & 4  & \leftarrow &  0    \\
              &   &    & \cdots &   &   &              
\end{array}
\end{array}
\]
}

\begin{table}[ht]
\caption{The Arnold numbers $v_{n,k}$.}
\begin{tabular}{c|cccccc|cccccc}
    \hline
$n$\textbackslash $k$
   & $-6$    & $-5$    & $-4$ & $-3$ & $-2$ & $-1$ & 1 & 2 & 3 & 4 &  5 &  6  \\
    \hline
          1   &   &  &  &    &    &  1  &  1  &    &    &   &   &   \\
          2   &   &  &  &    &  0 &  1  &  1  &  2 &    &   &   &   \\
          3   &   &  &  &  0 &  2 &  3  &  3  &  4 &  4 &   &  &    \\
          4   &   &  & 0 &  4 &  8 &  11 & 11  & 14 & 16 &  16 &  &  \\
          5   &   & 0 & 16 & 32 & 46 & 57 & 57  & 68 & 76 &  80 & 80 &  \\
          6   & 0 & 80 & 160 & 236 & 304 & 361 & 361 & 418 & 464 & 496 & 512 & 512   \\
      \hline    
\end{tabular}
\label{tab:Arnold-Springer-numbers}
\end{table}

A \emph{signed permutation} of $[n]$ is a bijection $\sigma$ of the set $[\pm n]:=\{-1,\dots,-n\}\cup\{1,\dots,n\}$ onto itself such that $\sigma(-i)=-\sigma(i)$ for any $i\in[\pm n]$. It will be denoted $\sigma=(\sigma_1,\dots,\sigma_n)$, where $\sigma_i=\sigma(i)$. A signed permutation $\sigma$ is a $\beta_n$-\emph{snake} if $\sigma_1>\sigma_2<\sigma_3>\cdots \sigma_n$ (cf. \cite[Section 8]{Arnold-1}, we change all values into their opposite). 
Let $B_n$ be the set of signed permutations of $[n]$, and let $\SSS_n\subset B_n$ be the set of $\beta_n$-snakes of size $n$. Note that $\#\SSS_n=2^nE_n$.

\medskip
\begin{thm} \label{thm:Arnold-Theorem}  {\rm (Arnold).} For $1\le |k|\le n$, we have 
\[
v_{n,k}:=\#\{\sigma\in\SSS_n: \sigma_1=k\}.
\]
\end{thm}

Following \cite{SZ}, a sequence of sets ($X_{n,k}$) is called an \emph{Arnold family} if $\#X_{n,k}=v_{n,k}$ for $1\le |k|\le n$.  
Josuat-Verg\`{e}s, Novelli, and Thibon \cite{JNT} studied $\beta_n$-snakes from the point of view of combinatorial Hopf algebras, and came up with an Arnold family in terms of valley-signed permutations. Recently, Shin and Zeng \cite{SZ} gave some Arnold families of objects: signed increasing 1-2 trees and a signed analogue of Andr\'{e} permutations of the second kind.

Hoffman \cite{Hoff-1} studied two sequences of polynomials $P_n$, $Q_n$ defined by
\begin{equation} \label{eqn:Hoffman-Pn-Qn}
\frac{d^n}{dx^n}\tan x =P_n(\tan x) \quad\mbox{and}\quad \frac{d^n}{dx^n}\sec x =Q_n(\tan x)\sec x,
\end{equation} 
for integer $n\ge 0$, and derived their exponential generating functions
\begin{equation} \label{eqn:EGF-for-Pn-Qn}
\sum_{n\ge 0} P_n(t)\frac{z^n}{n!}=\frac{\sin z+t\cdot\cos z}{\cos z-t\cdot\sin z} \quad\mbox{and}\quad
\sum_{n\ge 0} Q_n(t)\frac{z^n}{n!}=\frac{1}{\cos z-t\cdot\sin z}. 
\end{equation}
It was proved that $P_n(1)=2^nE_n$, $Q_n(1)=K(B_n)$, and $P_n(1)-Q_n(1)=K(D_n)$ \cite[Proposition 4.1]{Hoff-1}. Several of the initial polynomials are listed below.
\begin{align*}
P_1(t) &= 1+t^2,\\
P_2(t) &= 2t+2t^3, \\
P_3(t) &= 2+8t^2+6t^4, \\
P_4(t) &= 16t+40t^3+24t^5, \\
P_5(t) &= 16+136t^2+240t^4+120t^6. \\
& \\
Q_1(t) &= t,\\
Q_2(t) &= 1+2t^2, \\
Q_3(t) &= 5t+6t^3, \\
Q_4(t) &= 5+28t^2+24t^4, \\
Q_5(t) &= 61t+180t^3+120t^5.
\end{align*}

As a refinement of Arnold's result, we present double triangular arrays of polynomials, defined by recurrence, for the calculation of $P_n(t)$ and $Q_n(t)$. Namely, for $1\le |k|\le n$, we define the polynomials $V_{n,k}=V_{n,k}(t)$ by the following relations
\begin{equation} \label{eqn:Springer-poly-recurrence}
\begin{aligned}
V_{n,-k} &= V_{n,-k-1}+t^{-1} V_{n-1,k} \quad (n>k\ge 1),\\
V_{n,1} &= t^2\, V_{n,-1} \quad (n\ge 2),\\
V_{n,k}  &= V_{n,k-1}+t\, V_{n-1,-k+1} \quad (n\ge k>1),
\end{aligned}
\end{equation}
with $V_{1,1}=t^2$, $V_{1,-1}=1$, and $V_{n,-n}=0$ ($n\ge 2$). Several of the  polynomials $V_{n,k}(t)$ are listed in Table \ref{tab:Arnold-Springer-polynomials}. Note that $v_{n,k}=V_{n,k}(1)$.

\begin{table}[ht]
\caption{The Arnold--Hoffman polynomials $V_{n,k}(t)$.}
{\footnotesize
\begin{tabular}{c|ccccc}
    \hline
$n$\textbackslash $k$
        & 1 & 2 & 3 & 4 &  5  \\
    \hline
    1   &  $t^2$  &    &    &   &   \\
    2   &  $t^3$  &  $t+t^3$  &    &   &   \\
    3   &  $t^2+2t^4$  &  $2t^2+2t^4$ &  $2t^2+2t^4$ &   &  \\
    4   &  $5t^3+6t^5$  & $t+7t^3+6t^5$ & $2t+8t^3+6t^5$ &  $2t+8t^3+6t^5$ & \\
    5   &  $5t^2+28t^4+24t^6$  & $10t^2+34t^4+24t^6$ & $14t^2+38t^4+24t^6$ &  $16t^2+40t^4+24t^6$ & $16t^2+40t^4+24t^6$\\
    \hline
    \hline
$n$\textbackslash $k$
        & \multicolumn{1}{c}{$-5$} & \multicolumn{1}{c}{$-4$} & \multicolumn{1}{c}{$-3$} & \multicolumn{1}{c}{$-2$} &  \multicolumn{1}{c}{$-1$}      \\
    \hline
    1   &     &    &    &     &  1 \\
    2   &     &    &    &  0  &  $t$ \\
    3   &     &    &  0 &  $1+t^2$ & $1+2t^2$ \\
    4   &     & 0  & $2t+2t^3$ &  $4t+4t^3$ & $5t+6t^3$\\
    5   &  0  & $2+8t^2+6t^4$ & $4+16t^2+12t^4$ &  $5+23t^2+18t^4$ & $5+28t^2+24t^4$\\   
        \hline 
\end{tabular}
}
\label{tab:Arnold-Springer-polynomials}
\end{table}

One of our main results is a calculation of the polynomials $P_n(t)$, $Q_n(n)$ by the double triangular arrays ($V_{n,k}(t)$).

\begin{thm} \label{thm:Arnold-array-polynomials} For $1\le k\le n$, we have
\[
Q_n(t)=\frac{1}{t}\big(V_{n,1}+V_{n,2}+\cdots+V_{n,n} \big) \quad\mbox{and}\quad P_n(t)-t Q_n(t)=V_{n,-1}+V_{n,-2}+\cdots+V_{n,-n}.
\]
\end{thm}

Generalizing $Q_n$ in (\ref{eqn:Hoffman-Pn-Qn}), Josuat-Verg\`{e}s considered polynomials $Q^{(a)}_n$ such that
\[
\quad \frac{d^n}{dx^n}\sec^a x =Q^{(a)}_n(\tan x)\sec^a x.
\]
Of particular interest is the case $a=2$, denoted by $R_n=Q^{(2)}_n$. Since $\sec^2 x$ is the derivative of $\tan x$, it follows that $P_{n+1}(t)=(1+t^2)R_n(t)$ and $R_n(1)=2^nE_{n+1}$. Here are the initial polynomials of $R_n(t)$:
\begin{align*}
R_1(t) &= 2t,\\
R_2(t) &= 2+6t^2, \\
R_3(t) &= 16t+24t^3, \\
R_4(t) &= 16+120t^2+120t^4, \\
R_5(t) &= 272t+960t^3+720t^5.
\end{align*}

Josuat-Verg\`{e}s distinguished variants of $\beta_n$-snakes that realize the polynomials $P_n(t)$, $Q_n(t)$, and $R_n(t)$, interpreting the parameter $t$ by the statistic of sign changes \cite[Theorem 3.4]{JV}. 
The other aim of this paper is to provide some new Arnold families for $P_n(t)$, $Q_n(t)$ and related objects for $R_n(t)$, which are signed variants of  Andr\'{e} permutations and Simsun permutations, and find the statistics for the parameter $t$ in these families. In particular, we resolve an enumerative question raised by Shin and Zeng \cite[Conjecture 13]{SZ}.

\section{Preliminaries and Main results}
In this section we describe some Arnold families of objects that realize the polynomials $V_{n,k}$. 

\subsection{Complete increasing binary trees with empty leaves}
Josuat-Verg\`{e}s \cite[Section 4]{JV} derived the following objects from differential equations involving the series in (\ref{eqn:EGF-for-Pn-Qn}).

\begin{figure}[ht]
\begin{center}
\includegraphics[width=5.4in]{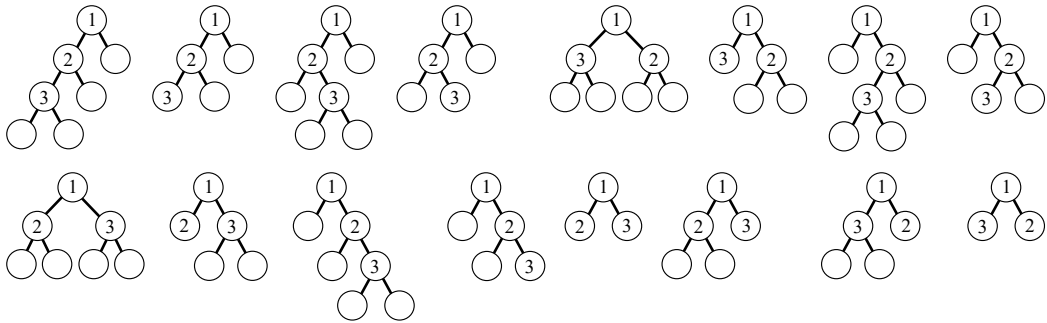}
\end{center}
\caption{\small The complete increasing binary trees with three labelled nodes.}
\label{fig:plane-labelled-binary-trees}
\end{figure}

A (plane) binary tree is \emph{complete} if every internal node has two children, i.e., a left child and a right child. 
Let $\T_n$ be the set of complete binary trees such that there are $n$ nodes labelled by integers from 1 to $n$ but some leaves can be non-labelled (these are called \emph{empty leaves}), and labels are increasing from the root to the leaves. Figure \ref{fig:plane-labelled-binary-trees} shows the sixteen complete increasing binary trees with three labelled nodes. 

For $\tau\in\T_n$, the \emph{rightmost path} of the tree $\tau$ is the unique sequence $(u_1,\dots,u_d)$ of vertices where $u_1$ is the root, $u_j$ is the right child of $u_{j-1}$ ($2\le j\le d$), and $u_d$ is a leaf.  We partition $\T_n$ into two subsets $\T^\ast_n$ and $\T^\circ_n$, where $\T^\ast_n$ ($\T^\circ_n$, respectively) consists of the trees whose rightmost leaf is labelled (non-labelled, respectively). Let $\emp(\tau)$ denote the number of empty leaves of $\tau$. Josuat-Verg\`{e}s obtained the following result \cite[Theorem 4.3]{JV}. We shall give a $q$-generalization of Theorem \ref{thm:Josuat-Verges-Theorem} in Section 6.

\begin{thm} \label{thm:Josuat-Verges-Theorem}  {\rm (Josuat-Verg\`{e}s).} For $n\ge 1$, we have
\[
P_n(t)=\sum_{\tau\in\T_n} t^{\emp(\tau)} \quad\mbox{and}\quad Q_n(t)=\sum_{\tau\in\T^\circ_n} t^{\emp(\tau)-1}.
\]
\end{thm}

Let $\ell(u)$ denote the label of a node $u\in\tau$. Let $\rmlab(\tau)$ denote the label of the last labelled node of the rightmost path of $\tau$, i.e., $\rmlab(\tau)=\ell(u_d)$ ($\ell(u_{d-1})$, respectively) if $u_d$ is labelled (non-labelled, respectively).  For $1\le k\le n$, let $\T^\ast_{n,k}\subset\T^\ast_n$ ($\T^\circ_{n,k}\subset\T^\circ_n$, respectively) be the subset of trees $\tau$ with $\rmlab(\tau)=k$. For example, as shown in Figure \ref{fig:plane-labelled-binary-trees}, $\#\T^\circ_{3,1}=4$, $\#\T^\circ_{3,2}=4$, and $\#\T^\circ_{3,3}=3$. Moreover, $\#\T^\ast_{3,3}=3$, $\#\T^\ast_{3,2}=2$, and $\#\T^\ast_{3,1}=0$. 

Notice that $\T^\ast_{1,1}$ consists of the tree with a single vertex (labelled 1), and that $\T^\circ_{1,1}$ consists of the tree with the root and two empty leaves. We present an interpretation of $V_{n,k}(t)$ in terms of the trees in $\T^\ast_{n,k}$ and $\T^\circ_{n,k}$. As a result, the sequence of sets $(\T^\ast_{n,1},\T^\ast_{n,2},\dots,\T^\ast_{n,n}$, $\T^\circ_{n,n},\T^\circ_{n,n-1},\dots,\T^\circ_{n,1})$ form an Arnold family, which provides an alternative proof of Arnold's result (see Corollary \ref{cor:beta-gamma}). 

\smallskip
\begin{thm} \label{thm:Arnold-family-plane-trees} For $1\le k\le n$, we have 
\[
V_{n,k}(t)=\sum_{\tau\in\T^\circ_{n,n-k+1}} t^{\emp(\tau)}\quad\mbox{and}\quad V_{n,-k}(t)=\sum_{\tau\in\T^\ast_{n,n-k+1}} t^{\emp(\tau)}.
\]
\end{thm}

\subsection{Signed Simsun permutations of type I}
For $\sigma\in\mathfrak{S}_n$, a \emph{descent} (\emph{ascent}, respectively) of $\sigma$ is a pair $(\sigma_i,\sigma_{i+1})$ such that $\sigma_i>\sigma_{i+1}$ ($\sigma_i<\sigma_{i+1}$, respectively) and $1\le i\le n-1$. A \emph{double descent} of $\sigma$ is a triple $(\sigma_{i},\sigma_{i+1},\sigma_{i+2})$ with $\sigma_{i}>\sigma_{i+1}>\sigma_{i+2}$ and $1\le i\le n-2$. 

\begin{defi} \label{def:Simsun-permutation} {\rm
A permutation $\sigma\in\mathfrak{S}_n$ is called a \emph{Simsun permutation} if for all $k$, the subword of $\sigma$ restricted to letters $\{1,\dots,k\}$ in the order they appear in $\sigma$ has no double descents.
}
\end{defi}

For $\sigma=(\sigma_1,\dots,\sigma_n)\in B_n$, let $|\sigma|:=(|\sigma_1|,\dots,|\sigma_n|)$. 
We distinguish some signed variants of Simsun permutations in Definitions \ref{def:Simsun-permutation-type-I-Rn} and \ref{def:Simsun-permutation-type-I}.

\begin{defi} \label{def:Simsun-permutation-type-I-Rn} {\rm
Let ${\RS}{\I}_n\subset B_n$ be the set of signed permutations $\sigma$ such that $|\sigma|$ is a Simsun permutation in $\mathfrak{S}_n$. These objects are called \emph{signed Simsun permutations of type I}. 
}
\end{defi}
Define the statistic $\npk$ of $\sigma\in {\RS}{\I}_n$ by
\begin{equation} \label{eqn:negative_peak}
\npk(\sigma):=\#\{\sigma_i:\sigma_i<0, |\sigma_i|>|\sigma_{i+1}|, 1\le i\le n-1\},
\end{equation}
i.e., the number of negative entries $\sigma_i$ such that $|\sigma_i|$ is a peak in $|\sigma|$. We obtain an interpretation of the polynomial $R_n(t)$.

\smallskip
\begin{thm} \label{thm:Rn(t)-interpretation} For $n\ge 1$, we have 
\[
R_n(t)=\sum_{\sigma\in{\RS}{\I}_n} t^{n-2\cdot\npk(\sigma)}. 
\]
\end{thm}

\begin{defi} \label{def:Simsun-permutation-type-I} {\rm
Let ${\RS}{\I}^{(B)}_n\subset{\RS}{\I}_n$ be the set of $\sigma$ that satisfies the following condition
\begin{equation} \label{eqn:condition-simsun-I-B}
 \sigma_i>0 \mbox{ for all $\sigma_i=\min\{|\sigma_i|,|\sigma_{i+1}|,\dots,|\sigma_n|\}$.}
\end{equation}
Let ${\RS}{\I}^{(D)}_n\subset{\RS}{\I}_n$ be the set of $\sigma$ that satisfies the following condition
\begin{equation} \label{eqn:condition-simsun-I-D}
 \sigma_n<0, |\sigma_n|>\sigma_{n-1}, \mbox{and $\sigma_i>0$ for all $\sigma_i=\min\{|\sigma_i|,|\sigma_{i+1}|,\dots,|\sigma_{n-1}|\}$.}
\end{equation}
}
\end{defi}
Here are these sets up to $n=3$.
\begin{align*}
{\RS}{\I}^{(B)}_1 &= \{(1)\};   \qquad{\RS}{\I}^{(B)}_2 = \{(1,2), (2,1), (-2,1) \}; \\
{\RS}{\I}^{(B)}_3 &= \{(1,2,3), (1,3,2), (1,-3,2), (2,1,3), (-2,1,3), (2,3,1), (-2,3,1), (2,-3,1),    \\
                  &\qquad (-2,-3,1), (3,1,2), (-3,1,2)\}. \\
{\RS}{\I}^{(D)}_1 &= \{(-1)\};   \qquad{\RS}{\I}^{(D)}_2 =  \{(1,-2)\}; \\
{\RS}{\I}^{(D)}_3 &= \{(1,2,-3), (2,1,-3), (-2,1,-3), (3,1,-2), (-3,1,-2)\}.                  
\end{align*}
For $1\le k\le n$, define
\[
{\RS}{\I}^{(B)}_{n,k}:=\{\sigma\in{\RS}{\I}^{(B)}_n: \sigma_n=k\} \quad\mbox{and}\quad {\RS}{\I}^{(D)}_{n,-k}:=\{\sigma\in{\RS}{\I}^{(D)}_n: \sigma_n=-k\}.
\]
Following Hetyei's signed Andr\'e permutations \cite{Hetyei}, the objects in ${\RS}{\I}^{(B)}_n$ were first studied by Shin and Zeng \cite{SZ}. They raised the following enumerative question (cf. \cite[Conjecture 13]{SZ}).

\begin{con} \label{con:conjecture-Shin-Zeng} {\rm (Shin--Zeng).} For $1\le k\le n$, we have
\[
v_{n,k}=\#{\RS}{\I}^{(B)}_{n,n-k+1}.
\]
\end{con}

We come up with the objects in ${\RS}{\I}^{(D)}_n$ and obtain the following realization of $V_{n,k}(t)$. It then follows that the sequence $({\RS}{\I}^{(D)}_{n,-1},\dots,{\RS}{\I}^{(D)}_{n,-n}, {\RS}{\I}^{(B)}_{n,n}, \dots,{\RS}{\I}^{(B)}_{n,1})$ is an Arnold family, which resolves the question in Conjecture \ref{con:conjecture-Shin-Zeng} (see Section 4).

\smallskip
\begin{thm} \label{thm:AF-simsun-I} For $1\le k\le n$, we have 
\[
V_{n,k}(t)=\sum_{\sigma\in{\RS}{\I}^{(B)}_{n,n-k+1}} t^{n+1-2\cdot\npk(\sigma)} \quad\mbox{and}\quad V_{n,-k}(t)=\sum_{\sigma\in{\RS}{\I}^{(D)}_{n,-n+k-1}} t^{n-1-2\cdot\npk(\sigma)}. 
\]
\end{thm}

\subsection{Signed Andr\'{e} permutations of type I}
Reformulated by Hetyei \cite[Definition 4]{Hetyei}, a permutation $\sigma\in\mathfrak{S}_n$ is an \emph{Andr\'{e} permutation of the second kind} if it is empty or satisfies the following conditions:
\begin{enumerate}
\item $\sigma$ has no double descent.
\item $(\sigma_{n-1},\sigma_n)$ is an ascent, i.e., $\sigma_{n-1}<\sigma_n$.
\item For $2\le i\le n-1$, if $\sigma_{i-1}>\sigma_{i}<\sigma_{i+1}$ then the minimum letter of $\omega_2$ is greater than the minimum letter of $\omega_4$ for the factorization $\sigma=\omega_1\,\omega_2\,\sigma_i\,\omega_4\,\omega_5$, where $\omega_2$, $\omega_4$ are maximal sequences of consecutive letters greater than $\sigma_i$.
\end{enumerate}
As noted by Shin and Zeng \cite{SZ}, the above definition is equivalent to the following definition.

\begin{defi} \label{def:Andre-permutation} {\rm
A permutation $\sigma\in\mathfrak{S}_n$ is called an \emph{Andr\'{e} permutation} if the subword of $\sigma$ restricted to letters $\{1,\dots,k\}$  has no double descent and ends with an ascent for all $1\le k\le n$.
}
\end{defi}

Following Hetyei \cite{Hetyei}, we define the following signed variants of Andr\'{e} permutations. 

\begin{defi} \label{def:signed-Andre-permutation-I} {\rm
Let $\AD{\I}_n\subset B_n$ be the set of signed permutations $\sigma$ containing the entry 1 such that $|\sigma|$ is an Andr\'{e} permutation in $\mathfrak{S}_n$. These objects are called \emph{signed Andr\'{e} permutations of type I}. Let $\AD{\I}^{(B)}_n\subset\AD{\I}_n$ be the set of $\sigma$ that satisfies the condition (\ref{eqn:condition-simsun-I-B}), and let $\AD{\I}^{(D)}_n\subset\AD{\I}_n$ be the set of $\sigma$ that satisfies the condition  (\ref{eqn:condition-simsun-I-D}). 
}
\end{defi}
The objects in $\AD{\I}^{(B)}_n$ were first defined by Hetyei \cite{Hetyei}. Here are these sets up to $n=4$.
\begin{align*}
{\AD}{\I}^{(B)}_2 &= \{(1,2)\};   \qquad{\AD}{\I}^{(B)}_3 = \{(1,2,3), (3,1,2), (-3,1,2) \}; \\
{\AD}{\I}^{(B)}_4 &= \{(1,2,3,4), (1,4,2,3), (1,-4,2,3), (3,1,2,4), (-3,1,2,4), (3,4,1,2),    \\
                  &\qquad (-3,4,1,2), (3,-4,1,2), (-3,-4,1,2), (4,1,2,3), (-4,1,2,3)\}. \\
{\AD}{\I}^{(D)}_2 &= \{(1,-2)\};   \qquad{\AD}{\I}^{(D)}_3 =  \{(1,2,-3)\}; \\
{\AD}{\I}^{(D)}_4 &= \{(1,2,3,-4), (3,1,2,-4), (-3,1,2,-4), (4,1,2,-3), (-4,1,2,-3)\}.                  
\end{align*}
For $2\le k\le n$, define
\[
\AD{\I}^{(B)}_{n,k}:=\{\sigma\in\AD{\I}^{(B)}_n: \sigma_n=k\} \quad\mbox{and}\quad \AD{\I}^{(D)}_{n,-k}:=\{\sigma\in\AD{\I}^{(D)}_n: \sigma_n=-k\}.
\]
Shin and Zeng \cite[Theorem 15]{SZ} established a bijection between $\AD{\I}^{(B)}_{n+1,k+1}$ and ${\RS}{\I}^{(B)}_{n,k}$, which can be extended to a bijection between $\AD{\I}^{(D)}_{n+1,-k-1}$ and ${\RS}{\I}^{(D)}_{n,-k}$ (see Theorem \ref{thm:bijection-AD-RS-I}).

\subsection{Signed Simsun permutations of type II} Let $\sigma\in B_n$. For $1\le k\le n$, let $\sigma_{[k]}$ denote the subword of $\sigma$ restricted to letters $\{\pm 1, \dots,\pm k\}$, i.e., obtained by removing the $n-k$ entries $\pm (k+1),\dots,\pm n$ from $\sigma$. 

\begin{defi} \label{def:Simsun-permutation-type-II-Rn} {\rm
Let ${\RS}{\II}_n\subset B_n$ be the set of signed permutations $\sigma$ such that the subword $\sigma_{[k]}$ contains no double descent for all $1\le k\le n$. These objects are called \emph{signed Simsun permutations of type II}. 
}
\end{defi}
Define the statistic $\nva$ of $\sigma\in {\RS}{\II}_n$ by
\begin{equation} \label{eqn:negative_valley}
\nva(\sigma):=\#\{\sigma_i:\sigma_i<0, \sigma_{i-1}>\sigma_i, |\sigma_{i-1}|<|\sigma_{i}|, 2\le i\le n\}.
\end{equation}
We obtain another interpretation of the polynomial $R_n(t)$.

\smallskip
\begin{thm} \label{thm:Rn(t)-alternative-interpretation} For $n\ge 1$, we have 
\[
R_n(t)=\sum_{\sigma\in{\RS}{\II}_n} t^{n-2\cdot\nva(\sigma)}. 
\]
\end{thm}

For $1\le k\le n$, the letter $k$ is called an \emph{augmenting element} of $\sigma$ if it is the last entry of the subword $\sigma_{[k]}$. Let $\gae(\sigma)$ denote the greatest augmenting element of $\sigma$. For example, the augmenting elements of $\sigma=(3,-1,2,4,7,5,8,-6)$ are $\{2,4,5\}$ and thus $\gae(\sigma)=5$.

\begin{defi} \label{def:Simsun-permutation-type-II} {\rm
Let ${\RS}{\II}^{(B)}_n\subset{\RS}{\II}_n$ be the set of $\sigma$ that satisfies the following condition
\begin{equation} \label{eqn:condition-simsun-II-B}
 \sigma_i>0 \mbox{ for all $\sigma_i=\min\{|\sigma_1|,|\sigma_2|,\dots,|\sigma_i|\}$.}
\end{equation}
Let ${\RS}{\II}^{(D)}_n\subset{\RS}{\II}_n$ be the set of $\sigma$ that satisfies the following condition
\begin{equation} \label{eqn:condition-simsun-II-D}
 \sigma_1<0, |\sigma_1|>\gae(\sigma), \mbox{and $\sigma_i>0$ for all $\sigma_i=\min\{|\sigma_2|,|\sigma_3|,\dots,|\sigma_i|\}$.}
\end{equation}
}
\end{defi}
Here are these sets up to $n=3$.
\begin{align*}
{\RS}{\II}^{(B)}_1 &= \{(1)\};   \qquad{\RS}{\II}^{(B)}_2 = \{(1,2), (1,-2), (2,1) \}; \\
{\RS}{\II}^{(B)}_3 &= \{(1,2,3), (1,-2,3), (1,2,-3), (1,3,2), (1,-3,2), (1,3,-2), (1,-3,-2), (2,1,3),     \\
                  &\qquad (2,3,1), (2,-3,1), (3,1,2)\}. \\
{\RS}{\II}^{(D)}_1 &= \{(-1)\};   \qquad{\RS}{\II}^{(D)}_2 = \{(-2,1)\}; \\
{\RS}{\II}^{(D)}_3 &= \{(-2,1,-3), (-2,3,1), (-3,1,2), (-3,1,-2), (-3,2,1)\}.                  
\end{align*}

We remark that the definition for ${\RS}{\II}^{(B)}_n$ is equivalent to the definition of signed Simsun permutations given by Ehrenborg and Readdy \cite{ER}: a signed permutation $\sigma\in B_n$ is a Simsun permutation if $\sigma_{[k]}$ contains no double descent and starts with a positive entry for all $1\le k\le n$. Ehrenborg and Readdy gave a recursive formula for the $cd$-index of the cubical lattice in terms of signed Simsun permutations \cite[Corollary 7.4]{ER}.

\smallskip
For $1\le k\le n$, define
\[
{\RS}{\II}^{(B)}_{n,k}:=\{\sigma\in{\RS}{\II}^{(B)}_n: \gae=k\} \quad\mbox{and}\quad {\RS}{\II}^{(D)}_{n,-k}:=\{\sigma\in{\RS}{\II}^{(D)}_n: \sigma_1=-k\}.
\]

We obtain another realization of $V_{n,t}(t)$, and hence the sequence $({\RS}{\II}^{(D)}_{n,-1},\dots,{\RS}{\II}^{(D)}_{n,-n}$, ${\RS}{\II}^{(B)}_{n,n},\dots,{\RS}{\II}^{(B)}_{n,1})$ is an Arnold family (see Section 5).

\smallskip
\begin{thm} \label{thm:AF-simsun-II} For $1\le k\le n$, we have
\[
V_{n,k}(t)=\sum_{\sigma\in{\RS}{\II}^{(B)}_{n,n-k+1}} t^{n+1-2\cdot\nva(\sigma)} \quad\mbox{and}\quad V_{n,-k}(t)=\sum_{\sigma\in{\RS}{\II}^{(D)}_{n,-n+k-1}} t^{n-1-2\cdot\nva(\sigma)}. 
\]
\end{thm}

\medskip
\subsection{Signed Andr\'{e} permutations of type II} We come up with a counterpart of $\RS{\II}_n$. 
\begin{defi} \label{def:signed-Andre-permutation-II} {\rm
Let $\AD{\II}_n\subset B_n$ be the set of signed permutations $\sigma$ containing the entry 1 such that the subword $\sigma_{[k]}$ contains no double descent and ends with an ascent for all $1\le k\le n$. These objects are called \emph{signed Andr\'{e} permutations of type II}. Let $\AD{\II}^{(B)}_n\subset\AD{\II}_n$ be the set of $\sigma$ that satisfies the following condition
\begin{equation} \label{eqn:condition-Andre-II-B}
 \sigma_i>0 \mbox{ for all $\sigma_i=\min\{|\sigma_1|,|\sigma_2|,\dots,|\sigma_i|\}$.}
\end{equation}
Let $\AD{\II}^{(D)}_n\subset\AD{\II}_n$ be the set of $\sigma$ that satisfies the following condition
\begin{equation} \label{eqn:condition-Andre-II-D}
 \sigma_1<0, |\sigma_1|>\sigma_n, \mbox{and $\sigma_i>0$ for all $\sigma_i=\min\{|\sigma_2|,|\sigma_3|,\dots,|\sigma_i|\}$.}
\end{equation}
} 
\end{defi}
Here are these sets up to $n=4$.
\begin{align*}
{\AD}{\II}^{(B)}_2 &= \{(1,2)\};   \qquad{\AD}{\II}^{(B)}_3 = \{(1,2,3), (1,-3,2), (3,1,2) \}; \\
{\AD}{\II}^{(B)}_4 &= \{(1,2,3,4), (1,2,-4,3), (1,-3,2,4), (1,4,2,3), (1,-4,2,3), (1,4,-3,2),   \\
                  &\qquad (1,-4,-3,2), (3,1,2,4), (3,4,1,2), (3,-4,1,2), (4,1,2,3)\}. \\
{\AD}{\II}^{(D)}_2 &= \{(-2,1,)\};   \qquad{\AD}{\II}^{(D)}_3 =  \{(-3,1,2)\}; \\
{\AD}{\II}^{(D)}_4 &= \{(-3,1,-4,2), (-3,4,1,2), (-4,1,2,3), (-4,1,-3,2), (-4,3,1,2)\}.                  
\end{align*}
For $2\le k\le n$, define
\[
{\AD}{\II}^{(B)}_{n,k}:=\{\sigma\in{\AD}{\II}^{(B)}_n: \sigma_n=k\} \quad\mbox{and}\quad {\AD}{\II}^{(D)}_{n,-k}:=\{\sigma\in{\AD}{\I}^{(D)}_n: \sigma_1=-k\}.
\]
We establish a bijection between ${\AD}{\II}_{n+1}$ and ${\RS}{\II}_n$, which induces bijections ${\AD}{\II}^{(B)}_{n+1,k+1}\rightarrow{\RS}{\II}^{(B)}_{n,k}$ and ${\AD}{\II}^{(D)}_{n+1,-k-1}\rightarrow{\RS}{\II}^{(D)}_{n,-k}$ (see Theorem \ref{thm:bijection-AD-RS-II}).

\section{Proof of Theorems \ref{thm:Arnold-array-polynomials} and \ref{thm:Arnold-family-plane-trees}}

For $1\le k\le n$, define the polynomials $T^\ast_{n,k}=T^\ast_{n,k}(t)$ and $T^\circ_{n,k}=T^\circ_{n,k}(t)$ by
\begin{equation}
T^\ast_{n,k}(t)=\sum_{\tau\in\T^\ast_{n,k}} t^{\emp(\tau)} \quad\mbox{and}\quad T^\circ_{n,k}(t)=\sum_{\tau\in\T^\circ_{n,k}} t^{\emp(\tau)}.
\end{equation}
In this section we shall prove that these polynomials $T^\ast_{n,k}$ and $T^\circ_{n,k}$ have the same recurrence relations of $(V_{n,k})$ as in (\ref{eqn:Springer-poly-recurrence}). Along with Josuat-Verg\`{e}s' bijection $\gamma:\T_n\rightarrow\SSS_n$, these recurrence relations provide an alternative proof of Arnold's result in Theorem \ref{thm:Arnold-Theorem}.

\subsection{Recurrence relations in terms of tree objects}

\medskip
\begin{pro}  \label{pro:recurrence-trees}  The following recurrence relations hold.
\begin{enumerate}
\item For $1< k\le n$,  $T^\ast_{n,k}=T^\ast_{n,k-1}+t^{-1}T^\circ_{n-1,k-1}$.
\item For $n\ge 2$,  $T^\circ_{n,n}=t^2\, T^\ast_{n,n}$,
\item For $1\le k< n$,  $T^\circ_{n,k}= T^\circ_{n,k+1}+t\,T^\ast_{n-1,k}$.
\end{enumerate}
\end{pro}

\begin{proof} 
(i) We shall establish a bijection $\psi$ between $\T^\ast_{n,k}$ and $\T^\ast_{n,k-1}\cup\T^\circ_{n-1,k-1}$. Given $\tau\in\T^\ast_{n,k}$, let $u,v$ be the nodes of $\tau$ with $\ell(u)=k-1$ and $\ell(v)=k$. Notice that $v$ is the rightmost leaf. The corresponding tree $\psi(\tau)$ is constructed as follows.
\begin{enumerate}
\item[(a)] If $v$ is not a child of $u$ then $\psi(\tau)$ is obtained from $\tau$ by interchanging the nodes $u$ and $v$. Notice that $\psi(\tau)\in\T^\ast_{n,k-1}$ and $\emp(\psi(\tau))=\emp(\tau)$.
\item[(b)] If $v$ is a child of $u$ then $\psi(\tau)$ is obtained from $\tau$ by replacing $v$ by an empty leaf and relabelling each node $x$ such that $\ell(x)\ge k+1$ by integer $\ell(x)-1$. Notice that $\psi(\tau)\in\T^\circ_{n-1,k-1}$ and $\emp(\psi(\tau))=\emp(\tau)+1$.
\end{enumerate}
The inverse map $\psi^{-1}$ can be established by the reverse operation.

(ii) There is an immediate bijection $\tau\mapsto \tau'$ of $\T^\ast_{n,n}$ onto $\T^\circ_{n,n}$ such that $\tau'$ is obtained from $\tau$ by adding two empty leaves to the rightmost leaf. Hence $\emp(\tau')=\emp(\tau)+2$.

(iii) We shall establish a bijection $\psi$ between $\T^\circ_{n,k}$ and $\T^\circ_{n,k+1}\cup\T^\ast_{n-1,k}$. Given $\tau\in\T^\circ_{n,k}$, let $v,w$ be the nodes of $\tau$ with $\ell(v)=k$ and $\ell(w)=k+1$. Notice that the right child of $v$ is an empty leaf. The corresponding tree $\psi(\tau)$ is constructed as follows.
\begin{enumerate}
\item[(a)] If $w$ is not a child of $v$ then $\psi(\tau)$ is obtained from $\tau$ by interchanging the nodes $v$ and $w$. Notice that $\psi(\tau)\in\T^\circ_{n,k+1}$ and $\emp(\psi(\tau))=\emp(\tau)$.
\item[(b)] If $w$ is a child of $v$ then it is the left child. There are two cases:
\begin{itemize}
\item $w$ is a leaf. Then $\psi(\tau)$ is obtained from $\tau$ by removing the two children of $v$ and relabelling each node $x$ such that $\ell(x)>k+1$ by integer $\ell(x)-1$. Notice that $\psi(\tau)\in\T^\ast_{n-1,k}$ and $\emp(\psi(\tau))=\emp(\tau)-1$.
\item $w$ is not a leaf. Let $\tau_1$ ($\tau_2$, respectively) be the left (right, respectively) subtree of $w$. Then $\psi(\tau)$ is obtained from $\tau$ by replacing the right child of $v$ by $w$, replacing the left subtree of $v$ ($w$, respectively) by $\tau_1$ ($\tau_2$, respectively), and letting the right child of $w$ be an empty leaf. Notice that $\psi(\tau)\in\T^\circ_{n,k+1}$ and $\emp(\psi(\tau))=\emp(\tau)$.
\end{itemize}
\end{enumerate} 
The inverse map $\psi^{-1}$ can be established by the reverse operation.
\end{proof}

\begin{figure}[ht]
\begin{center}
\includegraphics[width=1.6in]{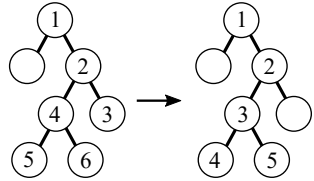}
\end{center}
\caption{\small An example for the map $\psi:\T^\ast_{n,k}\rightarrow\T^\circ_{n-1,k-1}$.}
\label{fig:map-(i)-binary-trees}
\end{figure}

\begin{exa} {\rm 
On the left in Figure \ref{fig:map-(i)-binary-trees} is a tree $\tau\in\T^\ast_{6,3}$ such that the rightmost leaf $v$, with $\ell(v)=3$, is a child of the node $u$ with $\ell(u)=2$. By case (b) of the map in the proof of Proposition \ref{pro:recurrence-trees}(i), the corresponding tree $\psi(\tau)\in\T^\circ_{5,2}$ is shown on the right.
}
\end{exa}

\begin{figure}[ht]
\begin{center}
\includegraphics[width=4in]{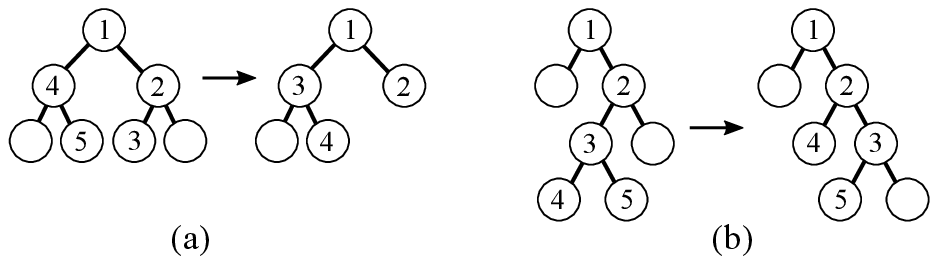}
\end{center}
\caption{\small An example for the map $\psi:\T^\circ_{n,k}\rightarrow\T^\circ_{n,k+1}\cup\T^\ast_{n-1,k}$.}
\label{fig:map-(iii)-binary-trees}
\end{figure}

\begin{exa} {\rm 
Figure \ref{fig:map-(iii)-binary-trees} shows two trees in $\T^\circ_{5,2}$ such that the node $w$ with $\ell(w)=3$ is a child of the node $v$ with $\ell(v)=2$. 
For the tree $\tau$ on the left in Figure \ref{fig:map-(iii)-binary-trees}(a) (\ref{fig:map-(iii)-binary-trees}(b), respectively), the node $w$ is a leaf (not a leaf, respectively). By case (b) of the map in the proof of Proposition \ref{pro:recurrence-trees}(iii), the corresponding tree $\psi(\tau)\in\T^\ast_{4,2}$ ($\T^\circ_{5,3}$, respectively) is shown on the right. 
}
\end{exa}

By Proposition \ref{pro:recurrence-trees} and the initial conditions $T^\circ_{1,1}=t^2$, $T^\ast_{1,1}=1$, and $T^\ast_{n,1}=0$ ($n\ge 2$),  we observe that the sequence of polynomials $(T^\ast_{n,1},T^\ast_{n,2},\dots,T^\ast_{n,n},T^\circ_{n,n},T^\circ_{n,n-1},\dots,T^\circ_{n,1})$ coincide with $(V_{n,-n},V_{n,-n+1},\dots,V_{n,-1},V_{n,1},V_{n,2},\dots,V_{n,n})$. The proof of Theorem \ref{thm:Arnold-family-plane-trees} is completed.

\medskip
\noindent
\emph{Proof of Theorem \ref{thm:Arnold-array-polynomials}.} By Theorems \ref{thm:Josuat-Verges-Theorem} and \ref{thm:Arnold-family-plane-trees}, we have

\begin{align*}
Q_n(t) &= \frac{1}{t}\sum_{\tau\in\T^\circ_n} t^{\emp(\tau)} \\
       &= \frac{1}{t}\sum_{k=1}^n\left(\sum_{\tau\in\T^\circ_{n,k}} t^{\emp(\tau)} \right)\\
       &= \frac{1}{t}\big(V_{n,1}+V_{n,2}+\cdots+V_{n,n} \big).
\end{align*}

\begin{align*}
P_n(t) &= \sum_{\tau\in\T_n} t^{\emp(\tau)} \\
       &= \sum_{\tau\in\T^\circ_n} t^{\emp(\tau)}+\sum_{\tau\in\T^\ast_n} t^{\emp(\tau)} \\
       &= tQ_n(t)+\big(V_{n,-1}+V_{n,-2}+\cdots+V_{n,-n} \big).
\end{align*}
\qed

The following relation is obtained by iterations of Proposition \ref{pro:recurrence-trees}(i).

\begin{cor} \label{cor:plane-trees-rightmost-leaf}
For $2\le k\le n$, we have $T^\ast_{n,k}=t^{-1}\big(T^\circ_{n-1,1}+T^\circ_{n-1,2}+\cdots+T^\circ_{n-1,k-1}\big)$.
\end{cor}

\subsection{Calculating $Q_n(t)$ by recurrence}
For $\tau\in\T_n$, the \emph{leftmost path} of the tree $\tau$ is the unique sequence $(u_1,\dots,u_d)$ of vertices where $u_1$ is the root, $u_j$ is the left child of $u_{j-1}$ ($2\le j\le d$), and $u_d$ is a leaf. Let $\T^{(L)}_n\subset\T_n$ denote the set of trees whose leftmost leaf is non-labelled. 

There is a bijection $\tau\mapsto\tau'$ of $\T^{(L)}_n$ onto $\T^\circ_n$ by flipping $\tau$  horizontally. It is easy to see that
\begin{equation} \label{eqn:(L)=circ}
\sum_{\tau\in\T^{(L)}_n} t^{\emp(\tau)}=\sum_{\tau\in\T^\circ_n} t^{\emp(\tau)}.
\end{equation}
For $1\le k\le n$, we define the following polynomials
\begin{equation} \label{eqn:(L)=circ+ast}
T^{(L)\ast}_{n,k}(t)=\sum_{\tau\in\T^{(L)}_n\cap\T^\ast_{n,k}} t^{\emp(\tau)} \quad\mbox{and}\quad T^{(L)\circ}_{n,k}(t)=\sum_{\tau\in\T^{(L)}_n\cap\T^\circ_{n,k}} t^{\emp(\tau)}.
\end{equation}
By the argument of the proof of Proposition \ref{pro:recurrence-trees}, we observe that these polynomials $T^{(L)\ast}_{n,k}(t)$ and $T^{(L)\circ}_{n,k}(t)$ have exactly the same recurrence relations as (i)-(iii) of Proposition \ref{pro:recurrence-trees}, along with the initial condition $T^{(L)\circ}_{n,1}=t^2$ and $T^{(L)\ast}_{n,1}=0$ ($n\ge 1$). By (\ref{eqn:(L)=circ}), (\ref{eqn:(L)=circ+ast}) and Theorem \ref{thm:Josuat-Verges-Theorem}, we have the following result for the calculation of $Q_n(t)$, which is analogous to the arrays for Anorld's $\gamma_n$-snakes (cf. \cite[Section 8]{Arnold-1}). Some of the initial polynomials of $T^{(L)\circ}_{n,k}(t)$ and $T^{(L)\ast}_{n,k}(t)$ are listed in Tables \ref{tab:Arnold-Springer-polynomials-Qn(t)-1} and \ref{tab:Arnold-Springer-polynomials-Qn(t)-2}, respectively.

\begin{cor} \label{cor:Qn(t)-by-T(L)_nk}
For $n\ge 1$, we have
\[
Q_n(t)=\frac{1}{t}\left(\sum_{k=1}^n T^{(L)\circ}_{n,k}(t)+ \sum_{k=1}^n  T^{(L)\ast}_{n,k}(t)\right).
\]
\end{cor}

\begin{table}[ht]
\caption{The polynomials $T^{(L)\circ}_{n,n-k+1}(t)$.}
{\footnotesize
\begin{tabular}{c|ccccc}
    \hline
$n$\textbackslash $k$
        & 1 & 2 & 3 & 4 &  5  \\
    \hline
    1   &  $t^2$  &    &    &   &   \\
    2   &  $t^3$  &  $t^3$  &    &   &   \\
    3   &  $2t^4$  &  $t^2+2t^4$ &  $t^2+2t^4$ &   &  \\
    4   &  $2t^3+6t^5$  & $4t^3+6t^5$ & $5t^3+6t^5$ &  $5t^3+6t^5$ & \\
    5   &  $16t^4+24t^6$  & $2t^2+22t^4+24t^6$ & $4t^2+26t^4+24t^6$ &  $5t^2+28t^4+24t^6$ & $5t^2+28t^4+24t^6$\\
            \hline 
\end{tabular}
}
\label{tab:Arnold-Springer-polynomials-Qn(t)-1}
\end{table}

\begin{table}[ht]
\caption{The polynomials $T^{(L)\ast}_{n,n-k+1}(t)$.}
{\footnotesize
\begin{tabular}{c|ccccc}
    \hline
$n$\textbackslash $k$
        & 5 & 4 & 3 & 2 &  1 \\
    \hline
    1   &     &    &    &     &  0 \\
    2   &     &    &    &  0  &  $t$ \\
    3   &     &    &  0 &  $t^2$ & $2t^2$ \\
    4   &     & 0  & $t+2t^3$ &  $2t+4t^3$ & $2t+6t^3$\\
    5   &  0  & $5t^2+6t^4$ & $10t^2+12t^4$ &  $14t^2+18t^4$ & $16t^2+24t^4$\\   
        \hline 
\end{tabular}
}
\label{tab:Arnold-Springer-polynomials-Qn(t)-2}
\end{table}

\subsection{An alternative proof of Arnold's results} Using Josuat-Verg\`{e}s' bijection $\gamma:\T_n\rightarrow\SSS_n$, we obtain Arnold's triangles about $\beta_n$-snakes and $\gamma_n$-snakes. We describe the bijection $\gamma$ below. 

Given $\tau\in\T_n$, read the word $\omega$ of the vertices of $\tau$ by an \emph{inorder} traversal. The word $\omega$ contains integers from 1 to $n$ and some letters, say $\circ$, for empty leaves.
Convert $\omega$ into $\omega'$ by replacing each element $i$ in $\omega$ by $n-i+1$ with sign $(-1)^{j+1}$, where $j$ is the number of letters $\circ$ after $i$, and removing all $\circ$'s. Then the corresponding $\beta_n$-snake $\gamma(\tau)$ is the reverse word of $\omega'$. 

\begin{figure}[ht]
\begin{center}
\psfrag{53c1c4c2c}[][][0.95]{$5,3,\circ,1,\circ,4,\circ,2,\circ$}
\psfrag{c154c23}[][][0.95]{$\circ,1,5,4,\circ,2,3$}
\includegraphics[width=3in]{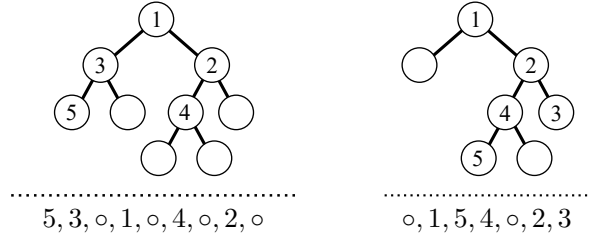}
\end{center}
\caption{\small The words read from the trees in Example \ref{exa:reading-words}.}
\label{fig:tree-to-snake}
\end{figure}

\smallskip
\begin{exa} \label{exa:reading-words} {\rm
For the tree $\tau\in\T^\circ_{5,2}$ shown on the left of Figure \ref{fig:tree-to-snake}, we have the word $\omega=(5,3,\circ,1,\circ,4,\circ,2,\circ)$. Converting $\omega$ into $\omega'=(-1,-3,5,-2,4)$, the corresponding $\beta_n$-snake is $\gamma(\tau)=(4,-2,5,-3,-1)$. 
Likewise, for the tree $\tau\in\T^\ast_{5,3}$ shown on the right of Figure \ref{fig:tree-to-snake}, we have $\omega=(\circ,1,5,4,\circ,2,3)$, $\omega'=(5,1,2,-4,-3)$, and hence $\gamma(\tau)=(-3,-4,2,1,5)$.
}
\end{exa}

Let $\SSS^\circ_n:=\{\sigma\in\SSS_n: (-1)^n\sigma_n<0\}$. The objects in $\SSS^\circ_n$ are called $\gamma_n$-snakes by Arnold \cite[Section 8]{Arnold-1}. For $1\le |k|\le n$, define
\[
\SSS^\circ_{n,k}=\{\sigma\in\SSS^\circ_n: \sigma_1=k\}.
\]

\begin{cor} \label{cor:beta-gamma} For $1\le k\le n$, we have
\begin{enumerate}
\item $\#\{\sigma\in\SSS_n:\sigma_1=k\}=\#\T^\circ_{n,n-k+1}$, and $\#\{\sigma\in\SSS_n:\sigma_1=-k\}=\#\T^\ast_{n,n-k+1}$.
\item $\#\SSS^\circ_{n,k}=\#\big(\T^{(L)}_n\cap\T^\circ_{n,n-k+1}\big)$, and $\#\SSS^\circ_{n,-k}=\#\big(\T^{(L)}_n\cap\T^\ast_{n,n-k+1}\big)$.
\end{enumerate}
\end{cor}

\begin{proof}
(i) For any $\sigma\in\SSS_n$ with $\sigma_1=k$ ($-k$, respectively), the map $\gamma^{-1}$ carries $\sigma$ to a tree $\gamma^{-1}(\sigma)$ in $\T^\circ_{n,n-k+1}$ ($\T^\ast_{n,n-k+1}$, respectively).

(ii) Given $\tau\in\T_n$, let $d$ be the number of labelled leaves of $\tau$. Then $\tau$ has $2(n-d)$ edges, i.e., $2(n-d)+1$ vertices. Hence the number of empty leaves of $\tau$ is $n+1-2d$. By the map $\gamma$, we observe that if the leftmost leaf of $\tau$ is non-labelled then the last entry of the $\beta_n$-snake $\gamma(\tau)$ is of sign $(-1)^{n+1}$. Hence $\gamma$ induces a bijection between $\T^{(L)}_n$ and $\SSS^\circ_n$. The assertion follows.
\end{proof}

Specializing the polynomials in Tables \ref{tab:Arnold-Springer-polynomials-Qn(t)-1} and \ref{tab:Arnold-Springer-polynomials-Qn(t)-2} at $t=1$, the initial values of $\#\SSS^\circ_{n,k}$ are shown in Table \ref{tab:Arnold-Springer-numbers-gamma}.

\begin{table}[ht]
\caption{Arnold's triangle ($\#\SSS^\circ_{n,k}$) for $\gamma_n$-snakes.}
\begin{tabular}{c|cccccc|cccccc}
    \hline
$n$\textbackslash $k$
   & $-6$    & $-5$    & $-4$ & $-3$ & $-2$ & $-1$ & 1 & 2 & 3 & 4 &  5 &  6  \\
    \hline
          1   &   &  &  &    &    &  0  &  1  &    &    &   &   &   \\
          2   &   &  &  &    &  0 &  1  &  1  &  1 &    &   &   &   \\
          3   &   &  &  &  0 &  1 &  2  &  2  &  3 &  3 &   &  &    \\
          4   &   &  & 0 & 3 &  6 &  8  &  8  & 10 & 11 & 11 &  &  \\
          5   &   & 0 & 11 & 22 & 32 & 40 & 40 & 48 & 54 &  57 & 57 &  \\
          6   & 0 & 57 & 114 & 168 & 216 & 256 & 256 & 296 & 328 & 350 & 361 & 361   \\
      \hline    
\end{tabular}
\label{tab:Arnold-Springer-numbers-gamma}
\end{table}

\section{Signed Simsun permutations of type I}
In this section we shall prove Theorem \ref{thm:Rn(t)-interpretation} and Theorem \ref{thm:AF-simsun-I} combinatorially using the following objects derived from complete increasing binary trees. 

\begin{figure}[ht]
\begin{center}
\includegraphics[width=4in]{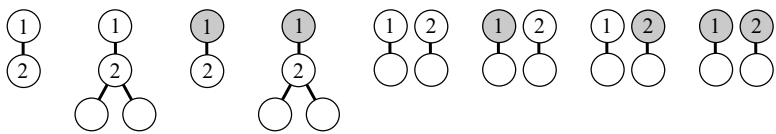}
\end{center}
\caption{\small The eight increasing forests with two labelled nodes.}
\label{fig:color-roots}
\end{figure}

\begin{defi} \label{def:plane-increasing-forest} {\rm
Let $\F_n$ be the set of rooted forests satisfying the following conditions.
\begin{enumerate}
\item Each root is colored either black or white. 
\item Each root has exactly one child, and each of the other internal nodes has exactly two ordered children.
\item There are $n$ nodes labelled by integers from 1 to $n$, but some leaves can be non-labelled (i.e., \emph{empty} leaves), and labels are increasing from each root down to the leaves.
\end{enumerate}
}
\end{defi}
Note that the order of the components of a forest is irrelevant. By convention, the components of a forest are arranged in increasing order of the roots.  Figure \ref{fig:color-roots} shows the set $\F_2$ of eight increasing forests with two labelled nodes.

Josuat-Verg\`{e}s (cf. \cite[Theorem 4.5]{JV}) proved that 
\begin{equation} \label{eqn:root-color-forest}
R_n(t)=\sum_{\pi\in\F_n} t^{\emp(\pi)}.
\end{equation}

\begin{figure}[ht]
\begin{center}
\includegraphics[width=3.6in]{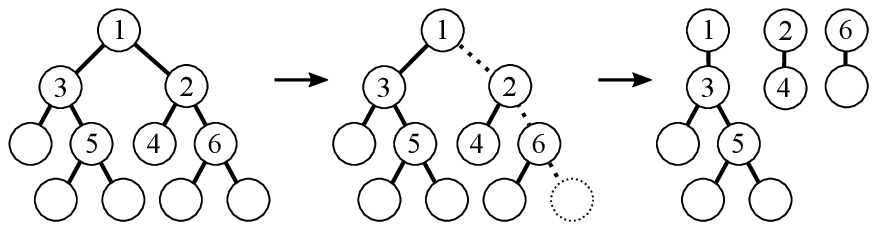}
\end{center}
\caption{\small The map $\mu:\T^\circ_{n,k}\rightarrow\F^{(W)}_{n,k}$.}
\label{fig:tree-to-forest}
\end{figure}

Let $\F^{(W)}_n\subset\F_n$ be the set of forests all of whose roots are white.
For $1\le k\le n$, let $\F^{(W)}_{n,k}\subset\F^{(W)}_n$ be the set of forests $\pi$ such that the root of the last component of $\pi$ is labelled with $k$. There is simple bijection $\mu:\T^\circ_{n,k}\rightarrow\F^{(W)}_{n,k}$: given $\tau\in\T^\circ_{n,k}$, remove the rightmost leaf and the edges in the path from the root to the rightmost leaf, and then the remaining components form the requested forest (see Figure \ref{fig:tree-to-forest}).

\subsection{A bijection $\varphi_1:\RS{\I}_n\rightarrow\F_n$}
To prove Theorem \ref{thm:Rn(t)-interpretation}, we shall establish the following bijection. 
  
\begin{pro} \label{pro:bijection-simsun-I-forest} There is a bijection $\varphi_1:\sigma\mapsto\pi$ of $\RS{\I}_n$ onto $\F_n$ with $\emp(\pi)=n-2\cdot\npk(\sigma)$. Moreover, for $1\le k\le n$, the map $\varphi_1$ induces a bijection between $\RS{\I}^{(B)}_{n,k}$ and $\F^{(W)}_{n,k}$.
\end{pro}

For $\sigma\in\RS{\I}_n$, an entry $\sigma_i$ (either positive or negative) is called a \emph{peak} (\emph{double-ascent element}, respectively) of $\sigma$ if $|\sigma_{i-1}|<|\sigma_i|>|\sigma_{i+1}|$ ($|\sigma_{i-1}|<|\sigma_i|<|\sigma_{i+1}|$, respectively), where we use the convention $\sigma_0=0$ and $\sigma_{n+1}=n+1$. This convention is also applied to every subword $\sigma_{[k]}$ of $\sigma$ ($1\le k\le n$). In particular, note that $\sigma_{[1]}=\pm 1$ is itself a double-ascent element.

To establish $\varphi_1$, we construct a sequence $\pi_1,\dots,\pi_n=\varphi_1(\sigma)$ of forests such that $\pi_j$ is determined by the subword $\sigma_{[j]}$ of $\sigma$.
For each $\pi_j$, the labelled leaves and the internal nodes with two empty leaves are called  \emph{terminal} nodes, and the internal nodes with exactly one empty leaf are called  \emph{intermediate} nodes. For example, let $\pi$ be the forest shown on the right of Figure \ref{fig:tree-to-forest}. The terminal (intermediate, respectively) nodes of $\pi$ are these nodes with labels $\{4,5\}$ ($\{3,6\}$, respectively). The construction of $\pi_1,\dots,\pi_n$ is described below.

\medskip
\noindent
{\bf Algorithm A}

The initial forest $\pi_1$ is the tree consisting of the root, labelled by integer 1, and one empty leaf, where the root is white (black, respectively) if $\sigma_{[1]}=1$  ($-1$, respectively). 
Suppose $\pi_{j-1}$ has been determined ($j\ge 2$), we locate the entry $x$ of the subword $\sigma_{[j]}$ with $|x|=j$, and then find out how $\sigma_{[j]}$ is obtained from $\sigma_{[j-1]}$:

\begin{enumerate}
\item  $x$ is inserted at the right of $\sigma_{[j-1]}$. Then $\pi_j$ is obtained from $\pi_{j-1}$ by adding a new component consisting of the root, labelled with $j$, and an empty leaf, where the root is white (black, respectively) if $x=j$ ($-j$, respectively).
\item  $x$ is inserted at the immediate left of a double-ascent element, say $y$, of $\sigma_{[j-1]}$. Find the intermediate node $v$ in $\pi_{j-1}$ with $\ell(v)=|y|$. Then $\pi_j$ is obtained from $\pi_{j-1}$ by labelling the empty leaf adjacent to $v$ by integer $j$. Moreover, the node labelled $j$ is assigned two empty leaves if $x>0$.
\item  $x$ is inserted at the immediate right of a peak, say $y$, of $\sigma_{[j-1]}$. Find the terminal node $v$ in $\pi_{j-1}$ with $\ell(v)=|y|$. If $y>0$ then $v$ has two empty leaves, and $\pi_j$ is obtained from $\pi_{j-1}$ by labelling the right child of $v$ by integer $j$. If $y<0$ then $v$ is a labelled leaf, and $\pi_j$ is obtained from $\pi_{j-1}$  by adding a node labelled $j$ (empty leaf, respectively) as the left (right, respectively) child of $v$. Moreover, the node labelled $j$ is assigned two empty leaves if $x>0$.
\end{enumerate}
Notice that the peaks (double-ascent elements, respectively) of $\sigma_{[j]}$ are in one-to-one correspondence with the terminal (intermediate, respectively) nodes of $\pi_{j}$. 

\begin{figure}[ht]
\begin{center}
\psfrag{F1}[][][0.85]{$\pi_1$}
\psfrag{F2}[][][0.85]{$\pi_2$}
\psfrag{F3}[][][0.85]{$\pi_3$}
\psfrag{F4}[][][0.85]{$\pi_4$}
\psfrag{F5}[][][0.85]{$\pi_5$}
\psfrag{F6}[][][0.85]{$\pi_6$}
\psfrag{F7}[][][0.85]{$\pi_7$}
\psfrag{F8}[][][0.85]{$\pi_8$}
\includegraphics[width=4.0in]{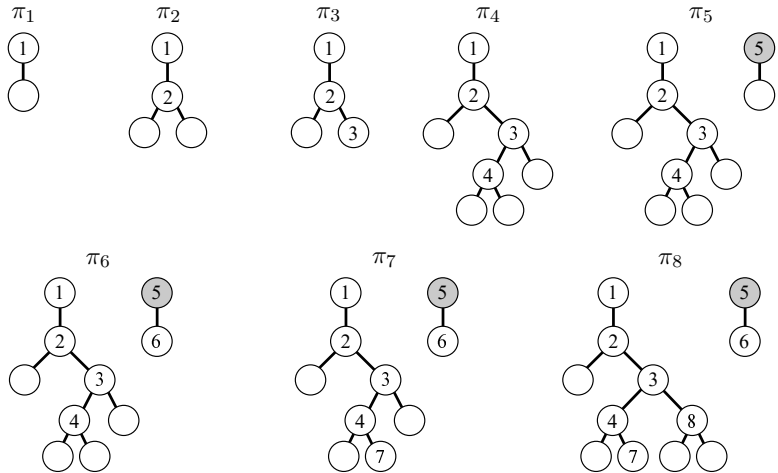}
\end{center}
\caption{\small The sequence of increasing forests for Example \ref{exa:simsun-to-forest}.}
\label{fig:simsun-forest-sequence-1}
\end{figure}

\begin{exa} \label{exa:simsun-to-forest} {\rm
Let $\sigma=(2,8,-3,4,-7,1,-6,-5)\in\RS{\I}_8$. The subwords $\sigma_{[j]}$ of $\sigma$ are shown in Table \ref{tab:construction-simsun-forest}, where the peaks and double-ascent elements of $|\sigma_{[j-1]}|$ are marked with $\ast$ and $\circ$, respectively. The sequence $\pi_1, \dots, \pi_8$ of forests for the construction of $\varphi_1(\sigma)$ are shown in Figure \ref{fig:simsun-forest-sequence-1}.  
}
\end{exa}

\begin{table}[ht]
\caption{The subwords $\sigma_{[j]}$ of $\sigma=(2,8,-3,4,-7,1,-6,-5)$.}
\centering
\begin{tabular}{rllll}
\hline
$x$ &  &\multicolumn{1}{c}{$|\sigma_{[j-1]}|$} & &\multicolumn{1}{c}{$\sigma_{[j]}$} \\
\hline \\ [-2ex]
1 & & & & $(1)$\\ [0.25ex]
2 & & $\begin{array}{rrrrrrrrr} 
 (0) & 1^{\circ} & (9) & & & & & & \end{array}$ & & $(2,1)$\\ [0.5ex]
$-3$ & & $\begin{array}{rrrrrrrrr} 
 (0) & 2^{\ast} & 1 & (9) & & & & & \end{array}$ & & $(2,-3,1)$\\ [0.5ex]
4 & & $\begin{array}{rrrrrrrrr} 
 (0) & 2^{\circ} & 3^{\ast} & 1 & (9) & & & & \end{array}$ & & $(2,-3,4,1)$\\ [0.5ex]
$-5$ & & $\begin{array}{rrrrrrrrr} 
 (0) & 2^{\circ} & 3^{\circ} & 4^{\ast} & 1 & (9) & &  \end{array}$ & & $(2,-3,4,1,-5)$\\ [0.5ex]
$-6$ & & $\begin{array}{rrrrrrrrr} 
 (0) & 2^{\circ} & 3^{\circ} & 4^{\ast} & 1 & 5^{\circ} & (9) &  & \end{array}$ & & $(2,-3,4,1,-6,-5)$\\ [0.5ex]
$-7$ & & $\begin{array}{rrrrrrrrr} 
  (0) & 2^{\circ} & 3^{\circ} & 4^{\ast} & 1 & 6^{\ast} & 5  & (9) & \end{array}$  & & $(2,-3,4,-7,1,-6,-5)$\\ [0.5ex]
8 & & $\begin{array}{rrrrrrrrr} 
  (0) & 2^{\circ} & 3^{\circ} & 4^{\circ} & 7^{\ast} & 1 & 6^{\ast} & 5  &  (9) \end{array}$  & & $(2,8,-3,4,-7,1,-6,-5)$\\ [0.5ex]
\hline
\end{tabular}
\label{tab:construction-simsun-forest}
\end{table}

To find $\varphi_1^{-1}$, given $\pi'\in\F_n$, we construct a sequence $\sigma'_{[1]},\dots,\sigma'_{[n]}=\varphi_1^{-1}(\pi')$ of words by the reverse operation of algorithm A.

\medskip
\noindent
{\bf Algorithm B1}
 
\begin{figure}[ht]
\begin{center}
\psfrag{+}[][][0.75]{$+$}
\psfrag{-}[][][0.75]{$-$}
\includegraphics[width=1.2in]{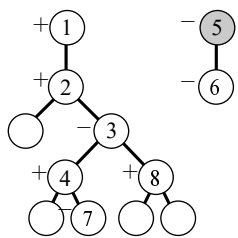}
\end{center}
\caption{\small The signs of the labelled nodes of an increasing forest.}
\label{fig:signed-labelled-forest}
\end{figure}

The first step is to assign a unique sign to each labelled node of $\pi'$ according to the following rules. 
By convention, each empty leaf is assumed to have a (virtual) value greater than $n$ so that if an empty leaf has a labelled sibling then the empty leaf is greater than its sibling. 
Define $\sign(u)=+1$ if the node $u$ satisfies one of the following conditions:
\begin{itemize}
\item $u$ is a white root, 
\item $u$ has two empty leaves, or
\item the left child of $u$ is greater than the right child, where the left child is possibly an empty leaf.
\end{itemize}
Moreover, define $\sign(u)=-1$ if the node $u$ satisfies one of the following conditions:
\begin{itemize}
\item $u$ is a black root, 
\item $u$ has no child, or
\item the left child of $u$ is less than the right child, where the right child is possibly an empty leaf.
\end{itemize}
For example, Figure \ref{fig:signed-labelled-forest} shows the signs of the labelled nodes of the forest $\pi_8$ in Figure \ref{fig:simsun-forest-sequence-1}.

Next, we construct a sequence $\pi'_1,\dots,\pi'_n$ of forests from $\pi'$ in reverse order. Let $\pi'_n=\pi'$. Suppose $\pi'_j$ has been determined, let $u$ be the node labelled by integer $j$ in $\pi'_j$. The forest $\pi'_{j-1}$ is obtained from $\pi'_j$ by replacing the node $u$ by an empty leaf and removing the children of $u$ if any. Moreover, (i) if $u$ is a root in $\pi'_j$ then remove the component of $u$, and (ii) if the parent of $u$ has sign $-1$ and two empty leaves in the resulting graph then remove these empty leaves. The construction of $\sigma'_{[1]},\dots,\sigma'_{[n]}$ is described below.

\medskip
\noindent
{\bf Algorithm B2}

By $\pi'_1$, we have the word $\sigma'_{[1]}=1$ ($-1$, respectively) if the root is white (black, respectively). Suppose $\sigma'_{[j-1]}$ has been determined ($j\ge 2$). 
Locate the node $u$ labelled $j$ in $\pi'_j$, with $\sign(u)\in\{1,-1\}$, and find out how $\pi'_j$ is obtained from $\pi'_{j-1}$:

\begin{enumerate}
\item  $u$ is the root of the newly added component. Then $\sigma'_{[j]}$ is formed by inserting the entry $\sign(u)j$ at the right of $\sigma'_{[j-1]}$.
\item  $u$ is a child of a node $v$, say with $\ell(v)=i$ and $\sign(v)\in\{1,-1\}$. There are two cases.
\begin{itemize}
\item $v$ is an intermediate node in $\pi'_{j-1}$. Then $\sigma'_{[j]}$ is formed by inserting the entry $\sign(u)j$ at the immediate left of the entry $\sign(v)i$ of $\sigma'_{[j-1]}$.
\item  $v$ is a terminal node in $\pi'_{j-1}$. Then $\sigma'_{[j]}$ is formed by inserting the entry $\sign(u)j$ at the immediate right of the entry $\sign(v)i$ of $\sigma'_{[j-1]}$.
\end{itemize}
\end{enumerate}

The bijection $\varphi_1:\RS{\I}_n\rightarrow\F_n$ is established.

\medskip
\noindent
\emph{Proof of Proposition \ref{pro:bijection-simsun-I-forest}.}
Given $\sigma\in\RS{\I}_n$, let $\pi=\varphi_1(\sigma)\in\F_n$ be the corresponding forest. Note that a labelled node $u$, say $\ell(u)=j$, is a leaf in $\pi$ whenever the entry $-j$ is a peak of $\sigma$. Hence the number of labelled leaves of $\pi$ equals $\npk(\sigma)$. Note that $\pi$ contains $2(n-\npk(\sigma))$ vertices. Hence $\emp(\pi)=n-2\cdot\npk(\sigma)$.

Moreover, an entry $x=\pm j$ of $\sigma$ is the last entry of the subword $\sigma_{[j]}$ if and only if $|x|$ is a right-to-left minimum of $|\sigma|$. 
For any $\sigma\in\RS{\I}^{(B)}_n$, by the condition (\ref{eqn:condition-simsun-I-B}) and algorithm A(i), we observe that all of the roots of $\varphi_1(\sigma)$ are white. Hence the map $\varphi_1$ induces a bijection between $\RS{\I}^{(B)}_{n,k}$ and $\F^{(W)}_{n,k}$. 
\qed

\smallskip
By (\ref{eqn:root-color-forest}) and Proposition \ref{pro:bijection-simsun-I-forest}, the proof of Theorem \ref{thm:Rn(t)-interpretation} is completed. Composing $\varphi_1$ with the bijection $\mu^{-1}:\F^{(W)}_{n,k}\rightarrow\T^\circ_{n,k}$ (in Figure \ref{fig:tree-to-forest}), we obtain the following result. 

\begin{cor} \label{cor:bijection-phi-I-B}
For $1\le k\le n$, there is a bijection $\phi^B_1:\sigma\mapsto\tau$ of $\RS{\I}^{(B)}_{n,k}$ onto $\T^\circ_{n,k}$ with $\emp(\tau)=n+1-2\cdot\npk(\sigma)$.
\end{cor}

By Theorem \ref{thm:Arnold-family-plane-trees} and Corollary \ref{cor:bijection-phi-I-B}, we have  $v_{n,k}=\#\T^\circ_{n,n-k+1}=\#\RS{\I}^{(B)}_{n,n-k+1}$, which proves Conjecture \ref{con:conjecture-Shin-Zeng}.

\subsection{A bijection $\phi^D_1:\RS{\I}^{(D)}_{n,-k}\rightarrow\T^\ast_{n,k}$}
On the basis of the previous bijection and the relation in Corollary \ref{cor:plane-trees-rightmost-leaf}, we obtain the following result. 

\begin{pro} \label{pro:bijection-simsun-I-forest-type-D}
For $2\le k\le n$, there is a bijection $\phi^D_1:\sigma\mapsto\tau$ of $\RS{\I}^{(D)}_{n,-k}$ onto $\T^\ast_{n,k}$ with $\emp(\tau)=n-1-2\cdot\npk(\sigma)$.
\end{pro}

\begin{proof}
Let $\sigma=(\sigma_1,\dots,\sigma_n)\in\RS{\I}^{(D)}_{n,-k}$. Note that $\sigma_n=-k$, and $\sigma_{n-1}=j$ for some $j<k$. We convert $\sigma$ into $\sigma'=(\sigma'_1,\dots,\sigma'_{n-1})\in\RS{\I}^{(B)}_{n-1,j}$ by setting
\[
\sigma'_i=\left\{\begin{array}{ll}
                    \sigma_i   &\mbox{if $|\sigma_i|<k$} \\              
                    \sigma_i-1 &\mbox{if $\sigma_i>k$} \\
                    \sigma_i+1 &\mbox{if $\sigma_i<-k$} 
                 \end{array}    
          \right.  
\]
for $1\le i\le n-1$. Notice that $\npk(\sigma')=\npk(\sigma)$. By Corollary \ref{cor:bijection-phi-I-B}, we construct the tree $\tau'=\phi^B_1(\sigma')\in\T^\circ_{n-1,j}$, where $\emp(\tau')=n-2\cdot\npk(\sigma')$. Then the corresponding tree $\phi^D_1(\sigma)$ is obtained from $\tau'$ by labelling the rightmost leaf by integer $k$, and relabelling each node $x$ such that $\ell(x)\ge k$ by integer $\ell(x)+1$. Hence $\emp(\phi^D_1(\sigma))=\emp(\tau')-1=n-1-2\cdot\npk(\sigma)$.

The inverse map of $\phi^D_1$ can be constructed by the reverse operation. 
\end{proof}

\begin{figure}[ht]
\begin{center}
\includegraphics[width=1.6in]{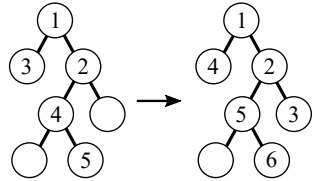}
\end{center}
\caption{\small The map $\psi^{-1}:\T^\circ_{n-1,j}\rightarrow\T^\ast_{n,k}$ for some $j<k$ in Example \ref{exa:type-D}.}
\label{fig:inverse-map-(i)}
\end{figure}

\begin{exa} \label{exa:type-D} {\rm 
Let $\sigma=(-4,1,5,-6,2,-3)\in\RS{\I}^{(D)}_{6,-3}$. Convert $\sigma$ into $\sigma'=(-3,1,4,-5,2)\in\RS{\I}^{(B)}_{5,2}$. By Corollary \ref{cor:bijection-phi-I-B}, we construct the tree $\phi^B_1(\sigma')\in\T^\circ_{5,2}$ shown on the left of Figure \ref{fig:inverse-map-(i)}. The corresponding tree $\phi^D_1(\sigma)\in\T^\ast_{6,3}$ is shown on the right.
}
\end{exa}

\smallskip
By Theorem \ref{thm:Arnold-family-plane-trees}, Corollary \ref{cor:bijection-phi-I-B} and Proposition \ref{pro:bijection-simsun-I-forest-type-D}, the proof of Theorem \ref{thm:AF-simsun-I} is completed.

\subsection{A bijection $\zeta_1:\AD{\I}_{n+1}\rightarrow\RS{\I}_n$.}  We describe a bijection between $\AD{\I}_{n+1}$ and $\RS{\I}_n$, which is extended from Shin and Zeng's \cite{SZ} bijection between Andr\'{e} permutation and Simsum permutations in $\mathfrak{S}_n$.

\begin{thm} \label{thm:bijection-AD-RS-I} There is a bijection $\zeta_1:\AD{\I}_{n+1}\rightarrow\RS{\I}_n$ such that for $1\le k\le n$, the map $\zeta_1$ induces bijections $\AD{\I}^{(B)}_{n+1,k+1}\rightarrow\RS{\I}^{(B)}_{n,k}$ and $\AD{\I}^{(D)}_{n+1,-k-1}\rightarrow\RS{\I}^{(D)}_{n,-k}$.
\end{thm} 

\begin{proof} Given $\sigma=(\sigma_1,\dots,\sigma_{n+1})\in\AD{\I}_{n+1}$, let $i_1<\cdots<i_d$ be the positions of the right-to-left minima of $|\sigma|$. Since $\sigma$ contains the entry 1 and $|\sigma|$ ends with an ascent, we have $\sigma_{i_1}=1$, $i_{d-1}=n$ and $i_d=n+1$. The corresponding word $\zeta_1(\sigma)=(\omega_1,\dots,\omega_n)\in\RS{\I}_n$ is defined by
\[
\omega_j = \left\{\begin{array}{ll}
                  \sigma_j-\dfrac{\sigma_j}{|\sigma_j|} &\mbox{if $j\not\in\{i_1,\dots,i_d\}$,} \\[3.2ex]
                  \sigma_{i_c}-\dfrac{\sigma_{i_c}}{|\sigma_{i_c}|} &\mbox{if $j=i_{c-1}$ for $c=2,\dots,d$.}
                  \end{array}
           \right.
\]
The inverse map $\zeta_1^{-1}$ can be defined reversely.

For $1\le k\le n$, note that $\omega_n=k$ ($-k$, respectively) if $\sigma_{n+1}=k+1$ ($-k-1$, respectively). Hence the assertions follows.
\end{proof}

\begin{exa} {\rm Let $\sigma=(3,6,1,2,9,-4,-8,5,-7)\in\AD{\I}_9$. The right-to-left minima of $|\sigma|$ are $(\sigma_3,\sigma_4,\sigma_6,\sigma_8,\sigma_9)=(1,2,-4,5,-7)$. Removing $1$ and moving the letters $ (2,-4,5,-7)$ to $(\sigma_3,\sigma_4,\sigma_6,\sigma_8)$ yields the word $(3,6,2,-4,9,5,-8,-7)$. Decrementing the absolute value of each entry by one, we have $\zeta_1(\sigma)=(2,5,1,-3,8,4,-7,-6)\in\RS{\I}_8$.
}
\end{exa}

\section{Signed Simsun permutations of type II}
In this section we shall prove Theorems \ref{thm:Rn(t)-alternative-interpretation} and \ref{thm:AF-simsun-II} by an approach similar to the proofs of Theorems \ref{thm:Rn(t)-interpretation} and \ref{thm:AF-simsun-I}. 

\subsection{A bijection $\varphi_2:\RS{\II}_n\rightarrow\F_n$}
We first establish the following bijection. 

\begin{pro} \label{pro:bijection-simsun-II-forest} There is a bijection $\varphi_2:\sigma\mapsto\pi$ of $\RS{\II}_n$ onto $\F_n$ with $\emp(\pi)=n-2\cdot\nva(\sigma)$.
Moreover, for $1\le k\le n$, the map $\varphi_2$ induces a bijection between $\RS{\II}^{(B)}_{n,k}$ and $\F^{(W)}_{n,k}$.
\end{pro}

For $\sigma\in\RS{\II}_n$, a descent pair $\sigma_i>\sigma_{i+1}$ is said to have a \emph{heavy top} (\emph{bottom}, respectively) if $|\sigma_i|>|\sigma_{i+1}|$ ($|\sigma_i|<|\sigma_{i+1}|$, respectively). Regarding the double-ascent elements of $\sigma$, we use the convention $\sigma_0=-(n+1)$ and $\sigma_{n+1}=n+1$, and apply to each subword $\sigma_{[k]}$ of $\sigma$.

To establish $\varphi_2$, we construct a sequence $\pi_1,\dots,\pi_n=\varphi_2(\sigma)$ of forests by the following procedure. For each $\pi_j$, we arrange the components with black (white, respectively) roots on the left (right, respectively) and in decreasing (increasing, respectively) order of the roots. By a \emph{singular} empty leaf of $\pi_j$ we mean a leaf that is the (only) child of a root, or a leaf whose sibling is labelled. 

\medskip
\noindent
{\bf Algorithm C}

The initial forest $\pi_1$ is the tree consisting of the root, labelled by integer 1, and an empty leaf, where the root is white (black, respectively) if $\sigma_{[1]}=1$ ($-1$, respectively). 
Suppose $\pi_{j-1}$ has been determined ($j\ge 2$), we locate the entry $x$ of the subword $\sigma_{[j]}$ with $|x|=j$, and then find out how $\sigma_{[j]}$ is obtained from $\sigma_{[j-1]}$:

\begin{enumerate}
\item  $x=j$ ($-j$, respectively) is inserted at the right (left, respectively) of $\sigma_{[j-1]}$. Then $\pi_j$ is obtained by adding a new component at the right (left, respectively) of $\pi_{j-1}$ consisting of the root and an empty leaf, where the root is white (black, respectively) and labelled by integer $j$.
\item  $x$ is inserted between an ascent pair, say $y<z$, of $\sigma_{[j-1]}$. 
There are two cases.
\begin{itemize}
\item $x=-j$. Note that by (i) $y\neq \sigma_0$ and that $y$ is a double-ascent element (otherwise there is a double descent in $\sigma_{[j]}$). Let $y$ be the $i$th double-ascent element of $\sigma_{[j-1]}$ from left to right for some $i$. Locate the $i$th singular empty leaf $u$ of $\pi_{j-1}$ from left to right. Then $\pi_j$ is obtained from $\pi_{j-1}$ by labelling $u$ by integer $j$.
\item $x=j$. Note that by (i) $z\neq \sigma_{n+1}$ and that $z$ is a double-ascent element. Let $z$ be the $i$th double-ascent element of $\sigma_{[j-1]}$ from left to right for some $i$. Locate the $i$th singular empty leaf $u$ of $\pi_{j-1}$ from left to right. Then $\pi_j$ is obtained from $\pi_{j-1}$ by labelling $u$ by integer $j$ and adding two empty leaves adjacent to $u$.
\end{itemize}
\item  $x$ is inserted between a descent pair, say $y>z$, of $\sigma_{[j-1]}$.  There are two cases.
\begin{itemize}
  \item If the descent pair has a heavy bottom (i.e., $|y|<|z|$) then $z<0$. Find the terminal node $v$ in $\pi_{j-1}$ with $\ell(v)=|z|$, and then add a node labelled $j$ (empty leaf, respectively) as the left (right, respectively) child of $v$.
  \item If the descent pair has a heavy top (i.e., $|y|>|z|$) then $y>0$. Find the terminal node $v$ in $\pi_{j-1}$ with $\ell(v)=y$, and then label the right child of $v$ by integer $j$.
\end{itemize}
In either case, the node labelled $j$ is assigned two empty leaves if $x>0$.
\end{enumerate}
Notice that the heavy elements (i.e., tops and bottoms) of $\sigma_{[j]}$ are in one-to-one correspondence with the terminal nodes of $\pi_{j}$, and that the number of double-ascent elements of $\sigma_{[j]}$ equals the number of singular empty leaves of $\pi_j$. 

\begin{figure}[ht]
\begin{center}
\psfrag{F1}[][][0.85]{$\pi_1$}
\psfrag{F2}[][][0.85]{$\pi_2$}
\psfrag{F3}[][][0.85]{$\pi_3$}
\psfrag{F4}[][][0.85]{$\pi_4$}
\psfrag{F5}[][][0.85]{$\pi_5$}
\psfrag{F6}[][][0.85]{$\pi_6$}
\psfrag{F7}[][][0.85]{$\pi_7$}
\psfrag{F8}[][][0.85]{$\pi_8$}
\includegraphics[width=4in]{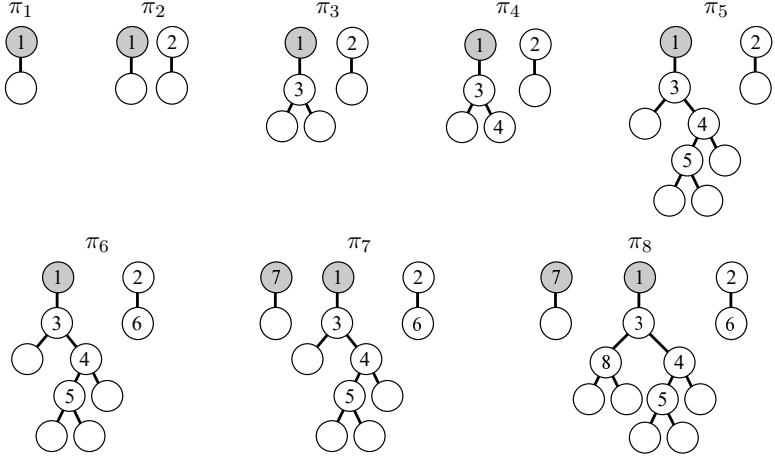}
\end{center}
\caption{\small The sequence of forests for Example \ref{exa:simsun-to-forest-II}.}
\label{fig:simsun-forest-sequence-2}
\end{figure}

\begin{exa} \label{exa:simsun-to-forest-II} {\rm
Let $\sigma=(-7,8,3,5,-4,-1,2,-6)\in\RS{\II}_8$. The subwords $\sigma_{[j]}$ of $\sigma$ are shown in Table \ref{tab:construction-simsun-II-forest}, where the heavy elements of $\sigma_{[j-1]}$ are marked with $\ast$, and the double-ascent elements are marked with $\circ$. The sequence $\pi_1, \dots, \pi_8$ of forests for the construction of $\varphi_2(\sigma)$ are shown in Figure \ref{fig:simsun-forest-sequence-2}. 
}
\end{exa}

\begin{table}[ht]
\caption{The subwords $\sigma_{[j]}$ of $\sigma=(-7,8,3,5,-4,-1,2,-6)$.}
\centering
\begin{tabular}{rllll}
\hline
$x$ &  &\multicolumn{1}{c}{$\sigma_{[j-1]}$} & &\multicolumn{1}{c}{$\sigma_{[j]}$} \\
\hline \\ [-2ex]
$-1$ & & & & $(-1)$\\ [0.25ex]
$2$ & & $\begin{array}{rrrrrrrrr} 
 (-9) & -1^{\circ} & (9) & & & & & & \end{array}$ & & $(-1,2)$\\ [0.5ex]
$3$ & & $\begin{array}{rrrrrrrrr} 
 (-9) & -1^{\circ} & 2^{\circ} & (9) & & & & & \end{array}$ & & $(3,-1,2)$\\ [0.5ex]
$-4$ & & $\begin{array}{rrrrrrrrr} 
 (-9) & 3^{\ast} & -1 & 2^{\circ} & (9) & & & & \end{array}$ & & $(3,-4,-1,2)$\\ [0.5ex]
$5$ & & $\begin{array}{rrrrrrrrr} 
 (-9) & 3 & -4^{\ast} & -1^\circ & 2^\circ & (9) & & & \end{array}$ & & $(3,5,-4,-1,2)$\\ [0.5ex]
$-6$ & & $\begin{array}{rrrrrrrrr} 
 (-9) & 3^\circ & 5^{\ast} & -4 & -1^\circ & 2^\circ & (9) & & \end{array}$ & & $(3,5,-4,-1,2,-6)$\\ [0.5ex]
$-7$ & & $\begin{array}{rrrrrrrrr} 
 (-9) & 3^\circ & 5^{\ast} & -4 & -1^\circ &  2  &  -6^\ast & (9) \end{array}$ & & $(-7,3,5,-4,-1,2,-6)$\\ [0.5ex]
$8$ & & $\begin{array}{rrrrrrrrr} 
 (-9) & -7^\circ & 3^\circ & 5^{\ast} & -4  & -1^\circ &  2  & -6^\ast & (9) \end{array}$ & & $(-7,8,3,5,-4,-1,2,-6)$\\ [0.5ex]
\hline
\end{tabular}
\label{tab:construction-simsun-II-forest}
\end{table}

To find $\varphi_2^{-1}$, given $\pi'\in\F_n$, we first assign a unique sign to each labelled node of $\pi'$ and construct a sequence $\pi'_1,\dots,\pi'_n$ of forests from $\pi'$ in reverse order, as we did in previous section by algorithm B1. Then we construct a sequence $\sigma'_{[1]},\dots,\sigma'_{[n]}=\varphi_2^{-1}(\pi')$ of words by the following procedure.

\medskip
\noindent
{\bf Algorithm D}

By $\pi'_1$, we have the initial word $\sigma'_{[1]}=1$ ($-1$, respectively) if the root is white (black, respectively). Suppose $\sigma'_{[j-1]}$ has been determined ($j\ge 2$). Locate the node $u$ labelled $j$ in $\pi'_j$, with $\sign(u)\in\{1,-1\}$, and then find out how $\pi'_j$ is obtained from $\pi'_{j-1}$:

\begin{enumerate}
\item  $u$ is the root of the newly added component. Then $\sigma'_{[j]}$ is formed by inserting the entry $j$ ($-j$, respectively) at the right (left, respectively) of $\sigma'_{[j-1]}$ if $u$ is white (black, respectively).
\item  $u$ is the $i$th singular empty leaf from left to right in $\pi'_{j-1}$. Then $\sigma'_{[j]}$ is formed by inserting the entry $j$ ($-j$, respectively) at the immediate left (right, respectively) of the $i$th double-ascent element of $\sigma'_{[j-1]}$ from left to right if $\sign(u)=1$ ($-1$, respectively).
\item  $u$ is a child of $v$ which is a terminal node in $\pi'_{j-1}$, say with $\ell(v)=y$ and $\sign(v)=1$ ($-1$, respectively). Then $\sigma'_{[j]}$ is formed by inserting the entry $\sign(u)j$ at the immediate right (left, respectively) of the entry $y$ ($-y$, respectively) of $\sigma'_{[j-1]}$.
\end{enumerate}

The bijection $\varphi_2:\RS{\II}_n\rightarrow\F_n$ is established. 

\smallskip
\noindent
\emph{Proof of Proposition \ref{pro:bijection-simsun-II-forest}.} 
For $\sigma\in\RS{\II}_n$, note that a node labelled $j$ is a leaf in $\varphi_2(\sigma)$ if and only if the entry $-j$ is a heavy bottom of $\sigma$. Hence the number of labelled leaves of $\varphi_2(\sigma)$ equals $\nva(\sigma)$. Note that $\varphi_2(\sigma)$ contains $2(n-\nva(\sigma))$ vertices. Hence $\emp(\sigma)=n-2\cdot\nva(\sigma)$.

Note that if $\varphi_2(\sigma)$ contains a black root, say labelled $j$, then the entry $-j$ is the first entry of the subword $\sigma_{[j]}$, and hence $j$ is a left-to-right minimum of $|\sigma|$. 
Restricted to $\sigma\in\RS{\II}^{(B)}_n$, by the condition (\ref{eqn:condition-simsun-II-B}) and algorithm C(i), we observe that all of the roots of $\varphi_2(\sigma)$ are white and that the last root is labelled by $\gae(\sigma)$. Hence the map $\varphi_2$ induces a bijection between $\RS{\II}^{(B)}_{n,k}$ and $\F^{(W)}_{n,k}$.
\qed

\smallskip
By (\ref{eqn:root-color-forest}) and Proposition \ref{pro:bijection-simsun-II-forest}, the proof of Theorem \ref{thm:Rn(t)-alternative-interpretation} is completed. 
Converting the forest $\varphi_2(\sigma)$ into a tree $\mu^{-1}(\varphi_2(\sigma))\in\T^\circ_n$, we obtain the following result.

\begin{cor} \label{cor:bijection-phi-II-B}
For $1\le k\le n$, there is a bijection $\phi^B_2:\sigma\mapsto\tau$ of $\RS{\II}^{(B)}_{n,k}$ onto $\T^\circ_{n,k}$ with $\emp(\tau)=n+1-2\cdot\nva(\sigma)$.
\end{cor}

\subsection{A bijection $\phi^D_2:\RS{\II}^{(D)}_{n,-k}\rightarrow\T^\ast_{n,k}$} On the basis of the previous bijection and the relation in Corollary \ref{cor:plane-trees-rightmost-leaf}, we obtain the following result.

\begin{pro} \label{pro:bijection-simsun-II-forest-type-D}
For $2\le k\le n$, there is a bijection $\phi^D_2:\sigma\mapsto\tau$ of $\RS{\II}^{(D)}_{n,-k}$ onto $\T^\ast_{n,k}$ with $\emp(\tau)=n-1-2\cdot\nva(\sigma)$.
\end{pro}

\begin{proof}
Let $\sigma=(\sigma_1,\dots,\sigma_n)\in\RS{\II}^{(D)}_{n,-k}$. Note that $\sigma_1=-k$, and $\gae(\sigma)=j$ for some $j<k$. We convert $\sigma$ into $\sigma'=(\sigma'_1,\dots,\sigma'_{n-1})\in\RS{\II}^{(B)}_{n-1,j}$ by setting
\[
\sigma'_i=\left\{\begin{array}{ll}
                    \sigma_{i+1}   &\mbox{if $|\sigma_{i+1}|<k$} \\
                    \sigma_{i+1}-1 &\mbox{if $\sigma_{i+1}>k$} \\
                    \sigma_{i+1}+1 &\mbox{if $\sigma_{i+1}<-k$} 
                 \end{array}    
          \right.  
\]
for $1\le i\le n-1$. Notice that $\nva(\sigma')=\nva(\sigma)$. By Corollary \ref{cor:bijection-phi-II-B}, we construct the tree $\tau'=\phi^B_2(\sigma')\in\T^\circ_{n-1,j}$, where $\emp(\tau')=n-2\cdot\nva(\sigma')$. Then the corresponding tree $\phi^D_2(\sigma)$ is obtained from $\tau'$ by labelling the rightmost leaf by integer $k$, and relabelling each node $x$ such that $\ell(x)\ge k$ by integer $\ell(x)+1$. Hence $\emp(\phi^D_2(\sigma))=\emp(\tau')-1=n-1-2\cdot\nva(\sigma)$.

The inverse map of  $\phi^D_2$ can be established by the reverse operation. 
\end{proof}

\smallskip
By Theorem \ref{thm:Arnold-family-plane-trees}, Corollary \ref{cor:bijection-phi-II-B} and Proposition \ref{pro:bijection-simsun-II-forest-type-D}, the proof of Theorem \ref{thm:AF-simsun-II} is completed. 

\subsection{A bijection $\zeta_2:\AD{\II}_{n+1}\rightarrow\RS{\II}_n$.}  We obtain an analogous result of Theorem \ref{thm:bijection-AD-RS-I}.

\begin{thm} \label{thm:bijection-AD-RS-II} There is a bijection $\zeta_2:\AD{\II}_{n+1}\rightarrow\RS{\II}_n$ such that for $1\le k\le n$, the map $\zeta_2$ induces bijections $\AD{\II}^{(B)}_{n+1,k+1}\rightarrow\RS{\II}^{(B)}_{n,k}$ and $\AD{\II}^{(D)}_{n+1,-k-1}\rightarrow\RS{\II}^{(D)}_{n,-k}$.
\end{thm} 

\begin{proof} Given $\sigma=(\sigma_1,\dots,\sigma_{n+1})\in\AD{\II}_{n+1}$, let $i_1<\cdots<i_d$ be the positions of augmenting elements of $\sigma$. Note that $\sigma_{i_1}=1$. Since every subword $\sigma_{[k]}$ of $\sigma$ ends with an ascent, $\sigma_{n+1}>0$ and $i_d=n+1$. Let $\omega=\zeta_2(\sigma)$. We construct the word $\omega=(\omega_1,\dots,\omega_n)$ as follows.

If $d=1$ then $\sigma_{n+1}=1$. The word $\omega$ is obtained from $\sigma$ by removing 1 and decrementing the absolute value of each remaining entry by one. If $d>1$, the word $\omega$ is defined by
\[
\omega_j = \left\{\begin{array}{ll}
                  \sigma_j-\dfrac{\sigma_j}{|\sigma_j|} &\mbox{if $j\not\in\{i_1,\dots,i_d\}$,} \\[2ex]
                  \sigma_{i_c}-1 &\mbox{if $j=i_{c-1}$ for $c=2,\dots,d$.}
                  \end{array}
           \right.
\]
We observe that any entry between $\sigma_{i_{c-1}}$ and $\sigma_{i_c}$ has its absolute value greater than $\sigma_{i_c}$ so that $\omega_{i_{c-1}}$ is an augmenting element in $\omega$. If not,  let $x$ be such an entry with the least absolute value, i.e, $|x|=\min\{|\sigma_{i_{c-1}+1}|,\dots,|\sigma_{i_c-1}|\}$ and $\sigma_{i_{c-1}}<|x|<\sigma_{i_c}$. Since $x$ is not an augmenting element, we have $x<0$ (otherwise, there is an entry on the right of $x$ with smaller absolute value, which is against the minimality of $|x|$), and it follows that the subword $\sigma_{[-x]}$ does not end with an ascent, a contradiction. 

We show that $\omega\in\RS{\II}_n$. Suppose $\omega_{[i]}$ has a double descent for some $1\le i\le n$. Let $\omega_a>\omega_b>\omega_c$ be three consecutive entries in $\omega_{[i]}$ for some $1\le a<b<c\le n$. At most one of $\omega_a,\omega_b,\omega_c$ is an augmenting element since the augmenting elements of $\omega$ are in increasing order.

If none of $\omega_a,\omega_b,\omega_c$ is augmenting then clearly $(\sigma_a,\sigma_b,\sigma_c)$ is a double descent in $\sigma_{[i+1]}$. Moreover, (i) if $\omega_c$ is augmenting then $\sigma_a=\omega_a+1$, $\sigma_b=\omega_b+1$ and $\sigma_c<\omega_c+1$; (ii) if $\omega_b$ is augmenting then $\omega_c<0$ and $|\omega_c|>\omega_b$, and it follows that $\sigma_a=\omega_a+1$, $0<\sigma_b<\omega_b+1$ and $\sigma_c=\omega_c-1<0$; (iii)  if $\omega_a$ is augmenting then $\omega_c<\omega_b<0$ and $|\omega_c|>|\omega_b|>\sigma_a$, and it follows that $0<\sigma_a<\omega_a+1$, $\sigma_b=\omega_b-1$ and $\sigma_c=\omega_c-1$. In each case $(\sigma_a,\sigma_b,\sigma_c)$ is a double descent in $\sigma_{[i+1]}$, which contradicts that $\sigma$ is a signed Andr\'{e} permutation of type II.

To find $\zeta_2^{-1}$, given $\omega'=(\omega'_1,\dots,\omega'_n)\in\RS{\II}_n$, let $\sigma'=\zeta_2^{-1}(\omega')$. If $\omega'$ contains no augmenting element then $\sigma'$ is obtained from $\omega'$ by incrementing the absolute value of each entry by one and adding the entry 1 to the right. Otherwise, let $i_1<\cdots<i_d$ be the positions of augmenting elements of $\omega'$. The word $\sigma'=(\sigma'_1,\dots,\sigma'_{n+1})$ is defined by
\[
\sigma'_j = \left\{\begin{array}{ll}
                  \omega'_j+\dfrac{\omega'_j}{|\omega'_j|} &\mbox{if $j\not\in\{i_1,\dots,i_d\}$,} \\[3ex]
                  1 &\mbox{if $j=i_1$,} \\[1ex]
                  \omega'_{i_{c-1}}+1 &\mbox{if $j=i_c$ for $c=2,\dots,d$,} \\[1ex]
                  \omega'_{i_d}+1 &\mbox{if $j=n+1$.}
                  \end{array}
           \right.
\]
Note that if $\sigma\in\AD{\II}^{(B)}_{n+1}$ with $\sigma_{n+1}=k+1$ then $\omega\in\RS{\II}_n$ with $\gae(\omega)=k$. Moreover, if $\sigma\in\AD{\II}^{(D)}_{n+1}$ with $\sigma_1=-k-1$ then $\omega\in\RS{\II}^{(D)}_n$ with $\omega_1=-k$. The assertions follow.
\end{proof}

\begin{exa} {\rm Let $\sigma=(4,-2,1,3,8,5,9,-7,6)\in\AD{\II}_9$. The augmenting elements of $\sigma$ are $(\sigma_3,\sigma_4,\sigma_6,\sigma_9)=(1,3,5,6)$. Removing $1$ and moving the letters $ (3,5,6)$ to $(\sigma_3,\sigma_4,\sigma_6)$ yields the word $(4,-2,3,5,8,6,9,-7)$. Decrementing the absolute value of each entry by one, we have $\zeta_2(\sigma)=(3,-1,2,4,7,5,8,-6)\in\RS{\II}_8$.
}
\end{exa}

\section{A $q$-generalization of Theorem \ref{thm:Josuat-Verges-Theorem}}

Let $D$ be the $q$-analogue of the differential operator acting on polynomials $f(t)$ by
\begin{equation}
(D\, f)(t) := \frac{f(qt)-f(t)}{(q-1)t}.
\end{equation}
Let $U$ be the operator acting on $f(t)$ by multiplication by $t$. Notice that the $q$-derivative $D(t^n)=[n]_q t^{n-1}$ and the communication relation $DU-qUD=1$ hold, where $[k]_q:=1+q+\dots +q^{k-1}$ for all positive integers $k$. Following \cite{JV}, we define the polynomials $P_n(q,t)$, $Q_n(q,t)$ and $R_n(q,t)$ algebraically by
\begin{equation}
P_n(q,t):=(D+UUD)^n t, \quad Q_n(q,t):=(D+UDU)^n 1, \quad R_n(q,t):= (D+DUU)^n 1.
\end{equation}
Here are several of the initial polynomials.
\begin{align*}
P_1(q,t) &= 1+t^2 \\
P_2(q,t) &= (1+q)t+(1+q)t^3 \\
P_3(q,t) &= 1+q+(2+3q+2q^2+q^3)t^2+(1+2q+2q^2+q^3)t^4, \\
&\\
Q_1(q,t) &= t \\
Q_2(q,t) &= 1+(1+q)t^2 \\
Q_3(q,t) &= (2+2q+q^2)t+(1+2q+2q^2+q^3)t^3, \\
  &\\
R_1(q,t) &= (1+q)t \\
R_2(q,t) &= 1+q+(1+2q+2q^2+q^3)t^2 \\
R_3(q,t) &= (2+5q+5q^2+3q^3+q^4)t+(1+3q+5q^2+6q^3+5q^4+3q^5+q^6)t^3. \\
\end{align*}

For $\tau\in\T_n$, let $\tau_0,\tau_1,\dots, \tau_n=\tau$ be the subtrees of $\tau$ such that $\tau_{j-1}$ is obtained from $\tau_j$ by replacing the node $u$ labelled $j$ by an empty leaf and removing the children of $u$ if any, and $\tau_0$ is an empty leaf. We associate $\tau$ with a sequence $(c_1,\dots,c_n)$ of integers, where $c_j$ is the number of empty leaves before the node $u$ (labelled $j$) in inorder traversal of $\tau_j$. Define the \emph{weight} of $\tau$ by $\w(\tau)=c_1+\cdots+c_n$.

\begin{figure}[ht]
\begin{center}
\psfrag{T1}[][][0.85]{$\tau_1$}
\psfrag{T2}[][][0.85]{$\tau_2$}
\psfrag{T3}[][][0.85]{$\tau_3$}
\psfrag{T4}[][][0.85]{$\tau_4$}
\psfrag{T5}[][][0.85]{$\tau_5$}
\includegraphics[width=4in]{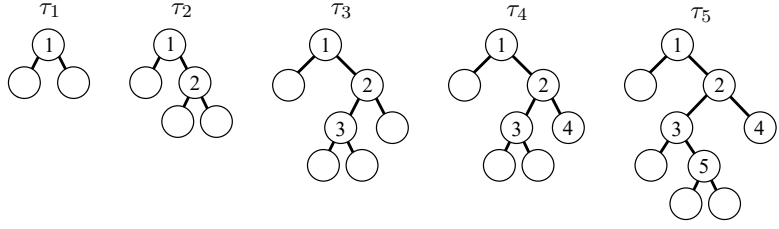}
\end{center}
\caption{\small The sequence of forests for Example \ref{exa:sub-trees}.}
\label{fig:subtree-sequence}
\end{figure}

\begin{exa}  \label{exa:sub-trees} {\rm
Figure \ref{fig:subtree-sequence} shows the subtrees $\tau_1,\dots,\tau_5$ of  a tree $\tau\in\T_5$, where $\tau=\tau_5$. The sequence associated with $\tau$ is $(c_1,\dots,c_5)=(0,1,1,3,2)$, and hence $\w(\tau)=7$.
}
\end{exa}  

\begin{thm} \label{thm:Pn(q,t)} For $n\ge 0$, we have
\[
P_n(q,t)=\sum_{\tau\in\T_n} q^{\w(\tau)}t^{\emp(\tau)}.
\]
\end{thm}

\begin{proof} Let $\tau_0$ be the tree consisting of an empty leaf, with weight $\w(\tau_0)=0$, and let $\T_0=\{\tau_0\}$. Note that $P_0(q,t)=t$, and hence the assertion holds for $n=0$.

By the proof of \cite[Theorem 4.3]{JV}, every $\tau\in\T_n$ can be established by constructing a sequence $\tau_0,\tau_1,\dots,\tau_n=\tau$ of trees, starting from $\tau_0$, such that $\tau_j$ is constructed from $\tau_{j-1}$ by labelling an empty leaf by integer $j$ and adding 0 or 2 empty leaves. Suppose $\tau_{j-1}$ has $k$ empty leaves for some $k\ge 0$, we observe that $\tau_j$ has a weight increment of $c_j$, i.e., $\w(\tau_j)=\w(\tau_{j-1})+c_j$, where $c_j$ is the number of empty leaves before the node labelled $j$ ($0\le c_j\le k-1$). Hence $\w(\tau)=c_1+\cdots+c_n$. Since $(D+UUD)t^k=[k]_q(t^{k-1}+t^{k+1})$, by induction we have
\[
P_n(q,t)=(D+UUD)P_{n-1}(q,t)=(D+UUD)\sum_{\tau\in\T_{n-1}} q^{\w(\tau)}t^{\emp(\tau)}=\sum_{\tau\in\T_n} q^{\w(\tau)}t^{\emp(\tau)}.
\]
The assertion follows.
\end{proof}

For $\pi\in\F_n$, the components of $\pi$ are arranged in increasing order of the roots. An order for the vertices of $\pi$ is given by from left to right traversing each component in inorder. Let $\pi_1,\dots, \pi_n=\pi$ be the subgraphs of $\pi$ such that $\pi_{j-1}$ is obtained from $\pi_j$ by replacing the node $u$ labelled $j$ by an empty leaf and removing the children of $u$ if any, and if $u$ is a root in $\pi_j$ then remove the component of $u$. The tree $\pi_1$ consists of a root, labelled 1, and an empty leaf. We associate $\pi$ with a sequence $(d_1,\dots,d_n)$ of integers, where $d_j$ is defined by
\[
d_j=\left\{\begin{array}{ll}
            c_j+1 &\mbox{if the node labelled $j$ is a black root in $\pi_j$} \\
            c_j &\mbox{otherwise,}
           \end{array}
   \right.          
\]
where $c_j$ is the number of empty leaves before the node labelled $j$. Define the weight of $\pi$ by $\w(\pi)=d_1+\cdots+d_n$. 

\begin{figure}[ht]
\begin{center}
\psfrag{F1}[][][0.85]{$\pi_1$}
\psfrag{F2}[][][0.85]{$\pi_2$}
\psfrag{F3}[][][0.85]{$\pi_3$}
\psfrag{F4}[][][0.85]{$\pi_4$}
\psfrag{F5}[][][0.85]{$\pi_5$}
\psfrag{F6}[][][0.85]{$\pi_6$}
\psfrag{F7}[][][0.85]{$\pi_7$}
\psfrag{F8}[][][0.85]{$\pi_8$}
\includegraphics[width=4.8in]{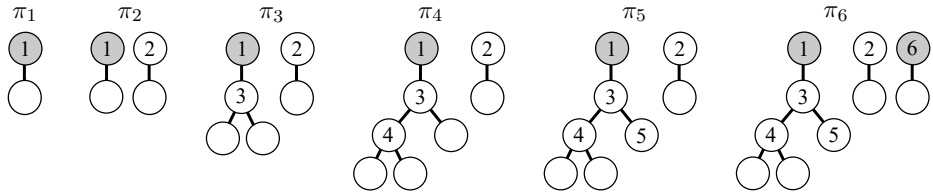}
\end{center}
\caption{\small The sequence of forests for Example \ref{exa:sub-forests}.}
\label{fig:subforest-sequence}
\end{figure}

\begin{exa} \label{exa:sub-forests} {\rm
Figure \ref{fig:subforest-sequence} shows the subgraphs $\pi_1,\dots,\pi_6$ of  a forest $\pi\in\F_8$, where $\pi=\pi_6$. The sequence associated with $\pi$ is $(d_1,\dots,d_6)=(1,1,0,0,2,4)$, and hence $\w(\pi)=8$.
}
\end{exa}   

\begin{thm} \label{thm:Rn(q,t)} For $n\ge 1$, we have
\[
R_n(q,t)=\sum_{\pi\in\F_n} q^{\w(\pi)}t^{\emp(\pi)},\qquad Q_n(q,t)=\sum_{\pi\in\F^{(W)}_n} q^{\w(\pi)}t^{\emp(\pi)}.
\]
\end{thm}

\begin{proof} Note that $\F_1$ consists of two trees, which contain a white (black, respectively) root and an empty leaf, with a weight of 0 (1, respectively).  Note that $R_1(q,t)=(D+DUU)1=(1+q)t$, and hence the first assertion holds for $n=1$.

For any $\pi\in\F_n$, we observe that $\pi$ can be established by constructing a sequence $\pi_1,\dots,\pi_n$ of forests, starting from $\pi_1\in\F_1$, such that $\pi_j$ is constructed from $\pi_{j-1}$ by either (i) labelling an empty leaf by integer $j$ and adding 0 or 2 empty leaves, or (ii) adding a new component consisting of the root, labelled $j$, and an empty leaf, where the root is either black or white. Suppose that $\pi_{j-1}$ has $k$ empty leaves for some $k\ge 0$, we observe that in the former case $\pi_j$ has a weight increment of $c_j$, where $c_j$ is the number of empty leaves before the node labelled $j$ ($0\le c_j\le k-1$), and in the latter case $\pi_j$ has a weight increment of $k$ ($k+1$, respectively) if the root labelled $j$ is white (black, respectively). Since $(D+DUU)t^k=[k]_qt^{k-1}+[k+2]_qt^{k+1}$, the first assertion is proved by the same induction argument as in the proof of Theorem \ref{thm:Pn(q,t)}. 

Moreover, if $\pi\in\F^{(W)}_n$ then in the latter case we have $\w(\pi_j)=\w(\pi_{j-1})+k$ since all roots of $\pi$ are white. Note that $Q_1(q,t)=(D+UDU)1=t$ and $(D+UDU)t^k=[k]_qt^{k-1}+[k+1]_qt^{k+1}$. The second assertion is proved similarly.

\end{proof}

\section{Conclusion}

In this paper, we make use of analogous recurrence relations of Arnold's method to calculate the polynomials $P_n(t)$, $Q_n(t)$ associated with successive derivatives of $\tan x$ and $\sec x$, established by Hoffman. With the benefit of the tree objects derived by Josuat-Verg\`es, our bijective proof in terms of the sets $\T^\ast_{n,k}$ and $\T^\circ_{n,k}$ provides a combinatorial insight for Arnold's enumeration of $\beta_n$-snakes. Using the tree objects, we describe some new Arnold's families in connection with Andr\'{e} permutations and Simsun permutations. 

Josuat-Verg\`{e}s, Novelli, and Thibon \cite{JNT} studied generalized snakes as alternating permutations with $r$ colors. The group of wreath product $\ZZ_r\wr \mathfrak{S}_n$ consists of all permutations $\omega$ of $[n]\times [0,r-1]$ such that $\omega(a,0)=(b,j)\Rightarrow \omega(a,i)=(b,i+j)$, where $i+j$ is computed modulo $r$, with composition of permutations as the group operation of $\ZZ_r\wr \mathfrak{S}_n$. 
Members $\omega=(\omega_1,\dots,\omega_n)$ of the group, called $r$-\emph{color permutations}, are represented by $(\sigma_1^{z_1},\cdots,\sigma_n^{z_n})$, where $(\sigma_1,\dots,\sigma_n)\in \mathfrak{S}_n$ and $(z_1,\dots,z_n)\in [0,r-1]^n$. Using the linear order
\[
1^{r-1}<\cdots<n^{r-1}<\cdots<1^{1}<\cdots<n^{1}<1<\cdots<n,
\]
an $r$-color alternating permutation $\omega$ is defined by the condition $\omega_1>\omega_2<\omega_3>\cdots\omega_n$. With $r=3$, the distribution of $r$-color alternating permutations with respect to the last entry is given in \cite[Proposition 8.1]{JNT}, up to $n=5$. We are interested in the corresponding tree object of $r$-color alternating permutations, which is helpful to study arithmetical properties of the distribution.

\end{document}